\font\steptwo=cmb10 scaled\magstep2

\magnification=\magstephalf
\settabs 18 \columns

\hsize=16truecm

\def\b{\bigskip}
\def\bb{\bigskip\bigskip}

\def\no{\noindent}
\def\r{\rightline}
\def\ce{\centerline}
\def\ve{\vfill\eject} 

\font\got=eufm10 scaled\magstep0
\font\bkbd=msbm10 scaled\magstep0

\def\g{\hbox{\got{g}}}
\def\hg{\hat{\g}}
\def\sl{\hbox{\got{sl}}}

\def\Cbb{\hbox{\bkbd{C}}}

\def\normal#1{\hbox{\rm :$#1$\rm :}}

\def\today{\ifcase\month\or January\or February\or March\or April\or 
May\or June\or July\or August\or September\or October\or November\or 
December\fi \space\number\day, \number\year }

\r \today
\b

\bb\bb 
{\ce {\steptwo  8-VERTEX  CORRELATION FUNCTIONS    }}\footnote{}{Mathematics Subject Classification (1991) 81R50, 70G50.}
\b
\ce {\steptwo AND} 
\b
\ce{\steptwo  TWIST COVARIANCE OF  q-KZ EQUATION}

\bb\bb
{\ce {C. Fr\o nsdal}}

{\ce {Physics Department, University of California, Los Angeles CA 90024, USA}}

\ce {and}

{\ce {A. Galindo}}
{\ce {Departamento de F\'{\i}sica Te\'orica, Universidad Complutense, 28040 Madrid, Spain}}

\bb\bb
 
 \no{\it ABSTRACT.}  We study the vertex operators $\Phi(z)$
associated with standard quantum groups.  The element $Z = RR^{\rm t}$
is a ``Casimir operator" for quantized Kac-Moody algebras and the
quantum Knizhnik-Zamolodchikov (q-KZ) equation is interpreted as the
statement $\normal{Z\Phi(z)} = \Phi(z)$.  We study the covariance of
the q-KZ equation under twisting, first within the category of Hopf
algebras, and then in the wider context of quasi Hopf algebras. We
obtain the intertwining operators associated with the elliptic
R-matrix and calculate the two-point correlation function for the
eight-vertex model.

\ve
\no{\bf {\steptwo   1. Introduction.}}
\b\no
In this paper we study the quantum Knizhnik-Zamolodchikov equation
[FR] for quasi Hopf algebras, with its covariance properties with
respect to twisting, and its relation to matrix elements of
intertwining operators. The conclusions bear on the interpretation of
the solutions of similar equations with exotic R-matrices. We
calculate the correlation functions for the 8-vertex model.

\b

\no {\it Correlation Functions for the Eight-Vertex Model.}

Baxter [Ba] introduced the trigonometric and elliptic quantum R-Matrix
for $\widehat {\sl(2)}$; this paper is mostly about the elliptic case,
and about the generalization [Be] to elliptic quantum $\widehat
{\sl(N)}$.  The trigonometric R-matrices found their interpretation in
terms of quantized Kac-Moody algebras, viewed as Hopf algebras; that
is, quantum groups [D1]. The elliptic R-matrices had, until recently,
not found their place in an algebraic framework. Surprisingly the
elliptic R-matrices also turned out to be related to quantized
Kac-Moody algebras, but with a quasi Hopf structure [Fr1,2]. More
precisely, the algebraic structure is the same as in the trigonometric
case, while the coproduct $\Delta$ of the trigonometric quantum group
is replaced by a new, deformed coproduct $\Delta_\epsilon$ (``elliptic
coproduct") that depends on a deformation parameter $\epsilon$. It can
be expressed as $\Delta_\epsilon = (F_\epsilon^{\rm t})^{-1} \Delta
F_\epsilon^{\rm t}$; the twistor $F_\epsilon$ must satisfy a cocycle
condition that has been solved to give an explicit expression for
$F_\epsilon$ as a power series in $\epsilon$. The quotient of the
elliptic quantum group, by the ideal generated by the center, is a
Hopf algebra; it is the quantization, in the sense of Drinfel'd, of
the classical, affine Lie bialgebra with elliptic r-matrix in the
classification of Belavin and Drinfeld [BD].

To understand the role of these elliptic quantum groups in the context
of integrable models and conformal field theory, we calculate the
correlation functions of the eight-vertex model.  The premise is that
Baxter's vertex operators can be interpreted mathematically as
intertwining operators for representations of quantized Kac-Moody
algebras [JM]; this is the interpretation that affords the most direct
link between statistical models and conformal field theory. Here we
define new intertwining operators in terms of the elliptic coproduct
and calculate the correlation functions that are associated with them;
that is, matrix elements of products of intertwining operators. We
find that these functions satisfy equations similar to the quantum
Knizhnik-Zamolodchikov equations of Frenkel and Reshetikhin [FR], but
that they can be described much more easily in terms of the familiar
correlation functions that govern the six-vertex model.

\b
\no {\it Twist Covariance.}

The larger issue is the question of the covariance of the q-KZ
equation under twisting in the category of quasi Hopf algebras. To
begin with, we point out that the q-KZ of Frenkel and Reshetikhin [FR]
can be easily generalized to all simple, affine quantum groups endowed
with what we call a ``standard" R-matrix: a universal R-matrix
(expressed as a series in Chevalley-Drinfeld generators, see
Definition 2.1.) that commutes with the Cartan subalgebra.
Reshetikhin [R] has described a highly specialized form of twisting
under which a standard R-matrix remains of standard type. From now on,
by the term ``twisting" we always have in mind a more radical twist
that transforms a standard R-matrix to a nonstandard or esoteric
R-matrix.

  A quantum group in the sense of this paper is a quantized, affine
Kac-Moody algebra $\hg$ based on a simple Lie algebra \g.  The
structure of coboundary Hopf algebra is given by a coproduct, an
antipode and a counit, but only the coproduct plays a direct role in
this paper. A coboundary Hopf algebra is a Hopf algebra $\hg$ with an
invertible element $R \in \hg\otimes\hg$ that satisfies the
Yang-Baxter relation and that intertwines the coproduct $\Delta$ with
its opposite $\Delta'$:
$$
R\Delta' = \Delta R.\eqno(1.1)
$$

The q-KZ equation is a holonomic system of difference equations that
 are satisfied by certain intertwining operators,
$$
\Phi, \Psi: V_{\mu,k} \rightarrow V(z) \otimes V_{\nu,k},
\eqno(1.2)$$
where $V_{\mu,k}$ and $V_{\nu,k}$ are irreducible, highest weight
$\hg$-modules of level $k$ and $V(z)$ is an evaluation module. The
intertwining property of $\Phi$ and of $\Psi$ is expressed as
$$
\Phi x = \Delta(x)\Phi,\quad \Psi x = \Delta'(x)\Psi,\eqno(1.3)
$$
for $x \in \hg$.  When $R$ is of standard type (Definition 2.1), then
the q-KZ equation for $\Psi$ takes the form
$$
(Z'-1)\Psi = 1,
\eqno(1.4)$$
where $Z'$ is a Casimir operator (acting in $V(z) \otimes V_{\nu,k}$)
  for $\hg$.  To define this operator let us express $R$ as
$$
R = R^i \otimes R_i,
$$ 
where we use the summation convention for the index $i$; then
formally,
$$
Z' = R^{\rm t}R, \quad R^{\rm t} := R_i \otimes R^i.\eqno(1.5)
$$
However, to make sense of an operator product such as $Z'\Psi$ it is
necessary to renormalize it.  The correct form of the q-KZ equation is
indeed (1.4), but with $Z'\Psi$ replaced by the normal-ordered product
$$
\normal{Z'\Psi} = R^{\rm t}(R^iq^{\hat H} \otimes 1)\Psi R_i,\eqno(1.6)
$$
where the factor $q^{\hat H}$ belongs to the Cartan subalgebra of
$\hg$.

 We study a deformation of the initial, standard quantum group,
implemented by twisting with an invertible element $F_\epsilon \in
\hg\otimes\hg$ that is a formal power series in a deformation
parameter $\epsilon$. The twisted quantities are
$$
\eqalign{&
R_\epsilon\ = (F_\epsilon^{\rm t})^{-1}RF_\epsilon,\quad
\Delta'_\epsilon = F_\epsilon\Delta' F_\epsilon^{-1},\cr &
\Psi_\epsilon = F_\epsilon^{-1}\Psi,\quad Z'_\epsilon =
F_\epsilon^{-1}Z'F_\epsilon,
\cr}
$$
and the twisted KZ equation is
$$
\normal{Z'_\epsilon\Psi_\epsilon} = \Psi_\epsilon;
$$
it has the same form as in the standard case. However, Eq.(1.6) is not
covariant; we mean by that it cannot be generalized by simply
replacing $R$ by $R_\epsilon $, since the expression
$$
R^{\rm t}_\epsilon(R^i_\epsilon \otimes 1)\Psi_\epsilon R_{\epsilon i}
$$
is not well defined. Instead, the correct expression for the
normal-ordered product is
$$
\normal{Z'_\epsilon\Psi_\epsilon} = F_\epsilon^{-1}\normal{Z'\Psi}
= F_\epsilon^{-1}R^{\rm t}(R^iq^{\hat H} \otimes 1)\Psi R_i.
$$
Therefore, though there is a clear sense in which ``the q-KZ equation"
is covariant, the normal-ordered product (1.6) is not.

This observation has analogous implications for correlation
function. To illustrate this, consider the two-point correlation
function $g(z_1,z_2) = \langle \Psi(z_1)\Psi(z_2) \rangle$.  In the
standard case the q-KZ equation reduces to
$$
g(q^{-k-g}z_1,z_2) = q^{A_1}R^{-1}(z_1,z_2)g(z_1,z_2).\eqno(1.7)
$$
The twisted correlation function obeys
$$
g_\epsilon(q^{-k-g}z_1,z_2) =
\bigl(F_\epsilon^{-1}(z_2,q^{-k-g}z_1)q^{A_1}R^{-1}(z_1,z_2)
F_\epsilon(z_2,z_1)\bigr)g_\epsilon(z_1,z_2),
$$
and this is {\it not} the same as Eq.(1.7) with $R$ replaced by
$R_\epsilon$.

This conclusion casts some light on the proposed generalization of of
the q-KZ equations for correlation functions. Integrability is assured
by the Yang-Baxter relation for the R-matrix.  It is natural to study
the equations that result from replacing the trigonometric R-matrix in
(1.7) and the rest, by more exotic R-matrices. Since this requires a
knowledge of such R-matrices in finite dimensional representations
only, it is possible, in particular, to use the elliptic R-matrix of
Baxter in this connection. As long as the elliptic quasi Hopf algebra
was not known, it was possible to speculate that the solutions of such
``elliptic q-KZ equations" relate in some way to (unknown) elliptic
intertwiners. Our conclusion is that this interpretation is not the
correct one.
 
\b
\no {\it Outline of the paper.}

Section 2 summarizes some facts about standard, universal R-matrices
and sets our notation.  Section 3 examines certain intertwining
operators and draws some conclusions (Proposition 3.1) that are used
later to determine the correct approach to regularizing operator
products.

Sections 4 and 5 present a view of the KZ and q-KZ equations. Both can
be interpreted very simply as eigenvalue equations, $\zeta \Phi = 0$
or $(Z-1)\Phi = 0$, for the Casimir operators $\zeta$ or $Z$ of affine
Kac-Moody or quantized, affine Kac-Moody algebras. Section 4 deals
with the classical KZ equation $\zeta \Phi = 0$; the effect of
different polarizations is discussed, as well as the invariance of the
operator $\zeta$ (Propositions 4.1 and 4.2). The quantum case is taken
up in Section 5; the correct normal-ordered action of the Casimir
elements $Z$ and $Z'$ on the intertwiners $\Phi$ and $\Psi$ is
established (Proposition 5.1), and the q-KZ equations are presented in
Eq.s (5.6) and (5.8).

Sections 6 and 7 explore the effect on intertwiners of twisting in the
categories of Hopf and quasi Hopf algebras.  In Section 6 we stress
the distinction between ``finite" and ``elliptic" twisting. The
twisted q-KZ equation is presented (Definition 6.3). In Section 7
quasi Hopf twisting is discussed and a recursion relation to actually
calculate the elliptic twistor is given.
 
Sections 8 and 9 apply the results to correlation functions.  In
Section 8 the classical and quantum q-KZ equations for correlation
functions are given; the effect of twisting is exhibited and a certain
lack of covariance is emphasized. In Section 9 the two-point
correlation function for the eight-vertex model is calculated, as well
as explicit expressions for the twisting matrix in the fundamental
representation of $\widehat {\sl(2)}$.

Finally, some auxiliary material is relegated to an Appendix.

\b
\no {\it Relation to other work.}

(1) Our original goal was to discover the enigmatic ``elliptic quantum
groups" and to use it to define and to calculate the correlation
functions for the eight-vertex model. This is precisely the
problematics of a series of paper by Jimbo, Miwa and others; see
especially the review [JM] and the paper [JMN].  These authors did not
have available the universal, elliptic R-matrix and did not anticipate
the fact that the algebraic structure of the elliptic quantum group
would turn out to be the same as in the trigonometric case. (Only the
coproduct is changed.)  They postulated a new algebraic structure, but
in the absence of a coproduct they could not define intertwiners. In
spite of this they did succeed in calculating correlation functions
that stand up to analysis and that reproduce some of Baxter's results
on the 8-vertex model.  Nevertheless, the correlation functions that
we here propose for the eight-vertex model are quite different.

(2) One of the most interesting aspects of the elliptic quantum group
is its quasi Hopf nature.  Quasi Hopf algebras, characterized by a
modified quantum Yang-Baxter relation, are basic to the
Knizhnik-Zamolodchikov-Bernard generalization of the KZ equation that
was discovered by Bernard [Ber]. This equation also arises in
connection with Felder's elliptic quantum groups [Fe].  However, these
developments are not concerned with highest weight matrix elements of
intertwiner operators, and the quasi Hopf algebras of Felder {\it et
al.} are not related to the elliptic R-matrices of Baxter and
Belavin. The new r-matrices discovered by Enriquez and Rubtsov [ER]
and by Frenkel, Reshetikhin and Semenov-Tian-Shansky [FRS] are of a
different sort. These interesting developments go beyond the
classification of classical r-matrices by Belavin and Drinfel'd [BD]
and are outside the scope of this paper.

\bb

\no {\bf {\steptwo 2. Standard, affine, universal,  quantum R-matrices.}}
\b\no
This section contains basic definitions and notation.

   Let $M,N$ be two finite sets,
$\varphi,\psi$  two maps,
$$
\eqalign{& \varphi : M\times M \rightarrow \Cbb, \cr & \psi
: M\times N \rightarrow \Cbb, \cr} \quad
\eqalign{a,b &\mapsto \varphi^{ab}, \cr a,\beta & \mapsto
H_a(\beta), \cr} 
$$
and $q$ a complex parameter.  Let ${\cal{A}}$ or
${\cal{A}}(\varphi,\psi)$ be the universal, associative, unital
algebra over {\Cbb} with generators $\{H_a\}_{ a\in M}, \{e_{\pm
\alpha}\}_{\alpha \in N}$, and relations
$$
\eqalignno{&[H_a,H_b]=0~, \quad  [H_a,e_{\pm\beta}] = \pm
H_a(\beta)e_{\pm\beta}, &   \cr
&[e_\alpha,e_{-\beta}]=\delta^\beta_\alpha
\bigl(q^{\varphi(\alpha,\cdot)}-q^{-\varphi(\cdot,\alpha)}\bigr), & 
  \cr}
$$
 with $\varphi(\alpha,\cdot)=\varphi^{ab}H_a(\alpha)H_b, 
\varphi(\cdot,\alpha)=\varphi^{ab}H_aH_b(\alpha)$ and $
q^{\varphi(\alpha, \cdot) + \varphi(\cdot,\alpha)} \neq 1, \alpha \in
N$.  The algebra of actual interest is a quotient ${\cal A}' = {\cal
A}/{\cal I}$ where ${\cal I}$ is a certain ideal; in this paper we
suppose that ${\cal I}$ is generated by a complete set of (quantized)
Serre relations among the $e_\alpha$'s and among the $e_{-\alpha}$'s;
then ${\cal A}'$ is a quantized (generalized) Kac-Moody algebra.  In
the case when ${\cal A}'$ is a quantized Kac-Moody algebra of affine
type, based on a simple Lie algebra \g, we sometimes write $\hg$ for
${\cal A}'$.  The ``Cartan subalgebra" ${\cal A}_0'$ is generated by
$\{H_a\}_{a \in M}$, extended by the inclusion of exponentials.

\b
\no {\bf Definition 2.1.} The standard, universal R-matrix has the form
$$ R = q^\varphi T = q^\varphi \sum_{n=0}^\infty t_n, \quad \varphi = 
\sum\varphi^{ab}  
H_a\otimes H_b, \eqno(2.1)
$$ where $t_0 = 1 \otimes 1, t_1 = \sum e_{-\alpha} \otimes e_\alpha$
(the sum is over the Serre generators, $\alpha \in N$) and $t_n$ has
the form
$$ t_n = t_{(\alpha)}^{(\alpha')}e_{-\alpha_1}...e_{-\alpha_n} \otimes
e_{\alpha_1'}...e_{\alpha_n'}.
\eqno(2.2)
$$ 
Sums over repeated indices are implied; the multi-index $(\alpha')$
 runs over the permutations of $(\alpha)$.
\b
The coefficients
$t_{(\alpha)}^{(\alpha')}
\in
\Cbb$ are essentially determined (the elements $t_n$ are determined 
uniquely) by the imposition of the Yang-Baxter relation,
$$
R_{12}R_{13}R_{23} = R_{23}R_{13}R_{12}.\eqno(2.3)
$$
It has been shown that, for a universal R-matrix of the type (2.1),
 this relation is equivalent to the recursion relation [Fr1]
$$ 
[  e_\gamma\otimes 1,t_n] = t_{n-1}(q^{\varphi(\gamma,.)}\otimes e_\gamma)
-(q^{-\varphi(.,\gamma)}\otimes e_\gamma)t_{n-1},
\eqno(2.4)
$$
 with the initial condition $t_0 = 1$. 
There is exactly one solution in ${\cal A}' \otimes {\cal A}'$.

We suppose now that ${\cal A}' = \hg$ is a quantized, affine Kac-Moody
algebra based on a simple Lie algebra \g. The coproduct is then
generated by the following formulas,
$$
 \Delta(e_\alpha) = 
 1 \otimes e_\alpha + e_\alpha
\otimes q^{\varphi(\alpha,.)},\quad \Delta(e_{-\alpha}) = 
 q^{-\varphi(.,\alpha)}  \otimes e_{-\alpha} + e_{-\alpha}
\otimes 1,\eqno(2.5)
$$
and $\Delta H_a = H_a \otimes 1 + 1 \otimes H_a$.  Let $\pi_1,\pi_2$
be finite dimensional representations of \g, and $\pi_i(z_i)$ the
associated evaluation representations of $\hg$ with spectral
parameters $z_i$. Let
$$
 R (z_1,z_2) := \pi_1(z_1) \otimes \pi_2(z_2)~ R. \eqno(2.6)
$$ 
The spectral parameters are regarded as formal variables; $R(z_1,z_2)$
is a formal power series $R_{12}(z_2/z_1)$ in $z_2/z_1$. The
effectiveness of the recursion relation (2.4) is illustrated in the
Appendix.

Finally, given $A = a^i \otimes a_i \in {\cal A}' \otimes {\cal A}'$,
we shall write $A^{\rm t} := a_i\otimes a^i$ and m$A := a^ia_i$.

\bb
\no{\bf {\steptwo 3.  Highest weight modules and intertwining operators.}}
\b\no
Let $V_\mu$ be an irreducible, finite dimensional, highest weight
$\hg$-module, and $V_{\mu,k} =
\bigoplus_{n \geq 0}V_{\mu,k}[-n]$ the associated
level $k$, highest weight, irreducible, graded $\hg$-module. The intertwining operators of 
greatest interest are imbeddings
$$
\Phi = \Phi(z): V_{\mu,k} \rightarrow V(z) \otimes V_{\nu,k},\eqno(3.1)
$$
where $V(z)$ is an evaluation module over $\hg$. The defining property
of $\Phi$ is
$$
\Phi x = \Delta (x) \Phi,
$$
for all $ x $ in $\hg$. 

We shall obtain some very essential information about the structure of
the intertwining operators.
\b

\no {\bf Proposition 3.1.} {\it Let $v$ be a homogeneous 
element of $V_{\mu,k}$. 
Then $\Phi v = \sum_n a_n \otimes b_n$, with $b_n \in V_{\nu,k}[-n]$,
where the sum is not, in general, finite.}
\b
\no {\it Proof.} It will be enough to verify that the sum is 
effectively infinite in one typical
case. Thus consider the quantized Kac-Moody algebra $\widehat{\sl(2)}$
with $V$ the fundamental representation.  In this case, for any $v \in
V_{\mu,k}$, $\Phi v$ takes the form
$$
\Phi v = \pmatrix{A(z)v\cr B(z)v\cr} = \pmatrix{A\cr B\cr} v.\eqno(3.2)
$$
Necessary  conditions to be satisfied by the operators 
$A,B: V_{\mu,k} \rightarrow V_{\nu,k}$ are
$$
\Delta(e_\alpha)\pmatrix{A\cr B\cr}v = 
\pmatrix{A\cr B\cr}e_\alpha v.\eqno(3.3)
$$
The first space is two-dimensional, with
$$
e_1 \otimes 1 = \kappa \pmatrix{0&1\cr0&0\cr},\quad e_0 \otimes 1 =
\kappa \pmatrix{0&0\cr z&0\cr}.
\eqno(3.4)
$$
The parameter $\kappa$ is related to $q, ~\kappa^2 = q - q^{-1}$.  In
full detail,
$$\eqalign{
&[e_0,A] = 0, \quad [e_1,A] = -\kappa q^{\varphi(1,.)}B,\cr
&[e_1,B] = 0,\quad [e_0,B] = -\kappa zq^{\varphi(0,.)}A.
\cr}\eqno(3.5)
$$
On the highest weight vector $v_0$ in $V_{\mu,k}$ we have
 $$
e_0Av_0 = 0, \quad e_1Bv_0=0,\quad e_0Bv_0 = -\kappa
zq^{\varphi(0,.)}Av_0,\quad e_1 A v_0 = -\kappa
q^{\varphi(1,.)}Bv_0,\eqno(3.6)
$$
 with a unique solution of the form 
$$
Av_0 = \sum_{n=0}^\infty z^nv_{2n}',\quad Bv_0 = \sum_{n = 0}^\infty
z^{n+1}v_{2n+1}',\eqno(3.7)
$$
with $v_0'$ a highest weight vector in $V_{\nu,k}$ and vectors $v_n'
\in V_{\nu,k}$ determined recursively by
$$
e_0v_{2n}' = e_1v_{2n+1}' = 0,\quad e_1v_{2n}' = -\kappa
q^{\varphi(1,.)}v_{2n-1}',\quad e_0v_{2n+1}' = -\kappa
q^{\varphi(0,.)}v_{2n}'.\eqno(3.8)
$$
The solutions have the form
$$
v'_{2n} = \sum_{\sigma \in S_{2n}}A_\sigma^{2n}
\,\sigma(e_{-1}e_{-0})^nv'_0,\quad v'_{2n+1} =
\sum_{\sigma \in S_{2n+1}}B_\sigma^{2n+1}
\sigma(e_{-1}e_{-0})^ne_{-0}v'_0, 
$$
where the sum is over all permutations of the generators. It is clear
that $v'_n \neq 0$ for all $n$ and the proposition is proved.
\b

We return to the general case, ${\cal A}'=\hg$ is a quantized
Kac-Moody algebra of affine type, based on a simple Lie algebra \g,
$V_{\mu,k}$ is a highest weight module over $\hg$, $V_i(z_i)$ finite
dimensional evaluation modules.
\b
\no {\it Remark 3.2.}
The product $\Phi_2\Phi_1$ is a compound map
$$
\Phi_1\quad \quad \quad \quad ~~\,\Phi_2\quad \,\,
$$
\vskip-0.9cm
$$
\Phi_2\Phi_1: ~ V_{\mu,k} \longrightarrow 
V_1 \otimes V_{\nu,k} \longrightarrow 
V_1 \otimes V_2 \otimes V_{\lambda,k}. \eqno(3.9)
$$
It has the property
$$
\Phi_2\Phi_1x = \Phi_2\Delta(x)\Phi_1 = 
({\rm id} \otimes\Delta)\Delta(x)\Phi_2\Phi_1.\eqno(3.10)
$$
By coassociativity of $\Delta$, $\Phi_2\Phi_1$ is an intertwiner of
the same type as $\Phi_1$ and $\Phi_2$:
$$
\Phi_2\Phi_1: V_{\mu,k} \rightarrow 
\bigl(V_1(z_1) \otimes V_2(z_2)\bigr) \otimes V_{\lambda,k}.
\eqno(3.11)
$$

Consequently, universal statements about intertwiners apply to
products of intertwiners as well.  This observation will be of use in
Section 8.  Of course, it does not apply in the quasi Hopf case
(Section 9.)

\ve

\no {\bf {\steptwo 4. The classical KZ equation.}}
\b\no
  The object
$$
Z = RR^{\rm t} \in {\cal A}' \otimes {\cal A}',\eqno(4.1)
$$
 if it exists, is invariant in the sense that it commutes with
$\Delta(x), \forall x\in {\cal A}'$. It plays the role of a Casimir
element for the quantized Kac-Moody algebra. Since the intertwiner
$\Phi$ projects on an irreducible representation, one expects that
there is $\langle Z \rangle \in \Cbb$ such that
$$
(Z - \langle Z \rangle)\Phi = 0.\eqno(4.2)
$$
We shall begin our study of this equation by considering its classical
limit. The result is Propositions (4.2) and (4.3). The important
concepts are normal ordering and ``polarization". Then we shall return
to the quantum case to show that (4.2) is the q-KZ equation of Frenkel
and Reshetikhin [FR]. (Section 5.)

The classical limit is defined by setting $q = {\rm e}^\eta$,
expanding in powers of $\eta$, and retaining the first nonvanishing
term. When ${\cal A}' = \hg$ is a quantized Kac-Moody algebra of
finite type, one finds that
$$
R = 1 + \eta r + O(\eta^2),\quad r = \varphi + 
\sum_{\alpha \in \Delta^+} E_{-\alpha} \otimes E_\alpha,\eqno(4.3)
$$
where the sum runs over the positive roots of \g. For simple roots one
has $ e_\alpha = \sqrt \eta (E_\alpha + O(\eta)); $ the others are
normalized so that the Casimir element in $\g\otimes\g$ takes the form
$$
{\cal C} = r + r^{\rm t}.
$$ 
 In the case of an untwisted affine loop algebra one gets
$$
R = 1 + \eta r + O(\eta^2),\quad r = \varphi + \sum_{\alpha \in
\Delta^+} E_{-\alpha}\otimes E_\alpha ~+ \sum_{n \geq 1}(z_2/z_1)^n
{\cal C},
\eqno(4.4)
$$
where $\Delta^+$ is the set of positive roots of the underlying Lie
algebra and where $z_1, z_2$ are the spectral parameters in the first,
resp. second space. It is important to keep in mind that this
expression is, until further development, nothing more than a formal
power series in $z_2/z_1$.  In terms of the basis
$$
E_{\pm \alpha}^n = z^nE_{\pm \alpha},\quad H_a^n = z^n H_a,\eqno(4.5)
$$
the expression for $r$ becomes
$$  
r = \varphi + E_{-\alpha} \otimes E_\alpha + \sum_{n\geq 1}{\cal C}^n
,\eqno(4.6)
$$
 with 
$$
{\cal C}^n =  K^{ab}H_a^{-n}
\otimes H_b^n + E_{-\alpha}^{-n} \otimes E_\alpha^n + E_\alpha^{-n} 
\otimes E_{-\alpha}^n,
\quad K^{ab} = (\varphi + \varphi^{\rm t})^{ab}.\eqno(4.7)
$$
Summation over $a,b$ and $\alpha \in \Delta^+$ will henceforth be
taken for granted.  Note that Eq.s (4.6) and (4.7) are valid in the
case of twisted loop algebras as well.

Returning to affine Kac-Moody algebras, it will be convenient to
change our conventions just a little.  Retain the above notation for
the loop algebra, so that, in particular,
$$
\varphi = \varphi^{ab} H_a \otimes H_b,\eqno(4.8)
$$
where the sum runs over the basis of the Cartan subalgebra of a simple
Lie algebra \g.  The form that characterizes the full, quantized
affine Kac-Moody algebra $\hg$ is
$$
\hat \varphi = \varphi + uc \otimes d + (1-u) d \otimes c,\eqno(4.9)
$$
where $d$ is the degree operator, $c$ is a basis for the central
extension and $u$ is a parameter.  For the full quantized Kac-Moody
algebra the limit is
$$
R = 1 + \eta \hat r + O(\eta^2), \quad \hat r = r + uc \otimes d+
(1-u) d \otimes c.\eqno(4.10)
$$

The classical limit of $Z$ is
$$
Z = 1 + \eta \zeta + O(\eta^2), \quad \zeta = \hat r + \hat r^{\rm
t}.\eqno(4.11)
$$
Formally,  
$$
\zeta =  \hat r + \hat r^{\rm t} = 
 \sum_{-\infty}^\infty {\cal C}^n + c\otimes d + d \otimes c.
$$
When both spaces are evaluation modules, where $c \mapsto 0$,
$$
\zeta = \sum_{n = -\infty}^{+\infty}(z_2/z_1)^n{\cal C}.\eqno(4.12)
$$
This sum becomes zero when projected on a quotient algebra of
meromorphic functions.
 
We  try to make sense out of the classical limit of (4.2), namely 
$$
(\zeta - \langle \zeta \rangle)\Phi = 0.
$$
By abuse of notation we retain the notation $\Phi$ for the classical
limit of the intertwiner.  Now the first space is an evaluation mode,
where $c$ vanishes, and if $c \mapsto k$ ($k$ is the level) on the
second space, then formally
 $$
\zeta\Phi(z) = kz{{\rm d}\over{\rm d}z}\Phi(z)  + 
 \sum_{n \geq 0} {\cal C}^{-n}\Phi(z) + \sum_{n > 0} {\cal
C}^n\Phi(z).\eqno(4.13)
$$
Let us introduce a uniform basis $ \{L_a\}$ for \g, so that the
 Casimir element takes the form ${\cal C} = L_a \otimes L_a$
 (summation implied). Then (in the untwisted case) (4.13) takes the
 form
$$
\zeta\Phi(z) = \biggl(kz{{\rm d}\over{\rm d}z}  + L_a    
\otimes\sum_{n \geq 0}z^nL_a^{-n}
 + L_a \otimes \sum_{n > 0}z^{-n}L_a^n\biggr)\Phi(z).\eqno(4.14)
$$
However, the significance of this formula is doubtful, as we shall
see.  This is the reason for the introduction of normal-ordered
products in [FR].
\b
\no {\it Polarization.}

It is usual, at this point of the development, to replace the operator
products by normal-ordered products.  It is a step that merits
comment. Normal-ordered operator products are introduced in field
theory when ordinary operator products fail to make sense. The typical
example is this product of destruction and creation operators:
$$
\bigl(\sum_n {\rm e}^{{\rm i}n\omega}a_n  
\bigr)\bigl(\sum_m {\rm e}^{-{\rm i}m\omega}a_m^*\bigr).
$$
When it is applied to the vacuum one gets
$$
\bigl(\sum_n {\rm e}^{{\rm i}n\omega}a_n  \bigr)\bigl(\sum_m {\rm e}^{-{\rm
i}m\omega}a_m^*\bigr)|0\rangle = \sum_{n = -\infty}^{+\infty}
|0\rangle,
$$
which is without meaning. The last term in (4.14) is of this kind; the
degree-decreasing operators in $\Phi$ correspond to the creation
operators, and the degree-increasing operators $L_a^n$ correspond to
the destruction operators. Collecting all terms of the same degree in
the product one gets a divergent series.  We want to avoid having to
interpret such infinite series, if it is possible.

Using the (classical) intertwining property of the  intertwiner,
$$
\sum_{n=1}^N z^{-n}( L_a \otimes L_a^n)\Phi = 
\sum_{n=1}^N z^{-n}(L_a \otimes 1)\Phi   L_a^n - 
(\sum_{n=1}^N L_a L_a \otimes 1)\Phi.\eqno(4.15)
$$
Passing with $N$ to infinity we encounter the meaningless expression
$\sum_{n > 0}L_aL_a$, an exact analogue of the divergent sum that is
thrown away when a field operator product is replaced by the
normal-ordered product.  It is tempting to redefine the operator
$\zeta$, by dropping this offensive term, thus
 $$
\zeta\Phi = kz{{\rm d}\over{\rm d}z}\Phi  +
  (L_a\otimes 1)\normal{J_a(z)\Phi}\quad (?),\eqno(4.16)
  $$
with
$$
J_a = J_a^+ + J_a^- = \sum_{n \geq 0}z^n L_a^{-n} + \sum_{n >
0}z^{-n}L_a^n,\eqno(4.17)
$$
and
$$
\normal{J_a\Phi} := J_a^+\Phi + \Phi J_a^-.\eqno(4.18) 
$$
This new operator is well defined on highest weight modules and it
will serve if it has the property that the formal expression (4.1) was
intended to assure; that is, if it is invariant. Actually it is,
almost.

 \b
\no{\bf Proposition 4.1.}   
{\it The covariant definition of the operator product
$\zeta\Phi(z)$ is
$$
\zeta\Phi(z) = (k+g)z{{\rm d}\over{\rm d}z}\Phi(z)  +
  (L_a\otimes 1)\normal{J_a(z)\Phi(z)},\eqno(4.19)
  $$
where $g$ is the dual Coxeter number of \g.}
\b
\no This result of [FR] is an analogue of Proposition 4.2 
that we prove below. The replacement of the factor $k$ by $k+g$, at
first sight somewhat mysterious, is thus required by covariance. For
$\sl(N)$, $g = N$.

\b
  We calculate the value $\langle \zeta \rangle$.  The operator $J_-$
 annihilates the highest weight vector $v_0$; therefore
$$
\zeta \Phi(z)v_0 = \bigl((k+g)z{{\rm d}\over{\rm d}z} + 
L_a \otimes J_a^+\bigr)\Phi(z)v_0.\eqno 
$$
In terms of the contravariant bilinear form (.,.), with $v'_0$ the
highest weight vector of $V_{\nu,k}$, one gets a $V(z)$-valued
function
$$ 
(v'_0,\,\Phi(z) v_0) =:\Phi_{v'_0v_0}(z) \in V(z),\eqno  
$$
and
$$\eqalign{
\bigl(v'_0,\,\zeta \Phi(z)v_0\bigr) &  = 
(k+g)z{{\rm d}\over{\rm d}z}\Phi_{v'_0v_0}(z) + (v'_0, L_a
\otimes J_a^+ ~\Phi(z)v_0).\cr}\eqno 
$$ 
  In $J_i^+$ only the zero mode contributes, and the second term
reduces to const.$\times\Phi_{v'_0v_0}(z)$. The constant has the value
$$
 {\textstyle{  1\over 2}}\bigl(C(\mu) - C(\nu) -  C(\pi)\bigr),\eqno(4.20) 
$$
where $C(\mu)$ is the value of the Casimir operator $C = {\rm
m}\,{\cal C}$ in $V_{\mu,k}[0]$. (Recall that if $A = a \otimes b \in
{\cal A}' \otimes {\cal A}'$, then m$A = ab \in {\cal A}'$.)
  
  We can reduce the value $\langle\zeta\rangle$ of $\zeta$  to zero
by choosing the grading of $\Phi$ according to 
$$
\Phi(z) = 
\sum_{n \in {\bf Z}}\Phi[n]z^{-n-(\mu|v\pi)}, \quad (\mu|\nu,\pi)  :=  
{1\over 2(k+g)}\bigl(C(\mu) - C(\nu) -  C(\pi)\bigr);
$$
 then for any weight vector $w$ in $V$,
$$
 \bigl(w \otimes v'_0, \Phi(z)v_0\bigr) = 
z^{-(\mu|v\pi)}\bigl(w\otimes v'_0,\Phi[0]v_0\bigr),\eqno(4.21)
$$
and
$$
\zeta\Phi(z) = (k+g)z{{\rm d}\over{\rm d}z}\Phi(z) 
+\normal{L_a \otimes J_a~\Phi(z)} =0.
\eqno(4.22) 
$$
This is the ``classical" Knizhnik-Zamolodchikov equation [KZ].
\b
\no {\it Alternative Polarizations.}

The polarization defined by (4.17), (4.18) is {\it ad hoc}. We have
the freedom of shifting any finite set of summands from $J^+$ to
$J^-$, as in
 $$
J_a = J_a^+ + J_a^- = \sum_{n > 0}z^n L_a^{-n} + \sum_{n \geq
0}z^{-n}L_a^n;
$$
the effect in this particular case is merely to change the sign of
$C(\pi)$ in (4.20). The result now agrees with [FR].

 Another polarization is suggested by (4.11),
$$
\zeta =  L_a \otimes J_a = 
L_a \otimes J_a^+ + L_a \otimes J^-_a = \hat r^{\rm t} + \hat r.
$$
Here we are dealing directly with the full Kac-Moody algebra,
including the $c,d$-terms in $\hat r$.  Formally, the intertwining
property gives
$$
\hat r\Phi = (\hat r^i \otimes \hat r_i)\Phi = 
-(\hat r^i\hat r_i \otimes 1)\Phi + 
(\hat r^i \otimes 1)\Phi\hat r_i,
$$
and
$$
\hat r^i\hat r_i = {\textstyle {1\over 2}}\bigl(\sum   
C^n  + c \otimes d + d \otimes c\bigr)+ 
{\textstyle {1\over 2}}[\hat r^i,\hat r_i].\eqno(4.23)
$$
The first term on the right hand side of this last equation, though
meaningless, looks like it may be a scalar, and thus
ignorable. Proceeding heuristically up to Proposition 4.2, we begin by
dropping this term.  The other term is an element $\hat H$ of the
Cartan subalgebra of $\hg$,
$$
\hat H = {\textstyle {1\over 2}}[\hat r_i,\hat r^i];\eqno(4.24) 
$$
it is determined up to an additive central element by
$$
[\hat H,e_\alpha] = [e_\alpha,\hat r^i\,\hat r_i] = {\rm{\rm
m}}[e_\alpha \otimes 1 + 1 \otimes e_\alpha, \hat r] = {\rm
m}\bigl(\varphi(\alpha,.)
\wedge e_\alpha\bigr) =  \varphi(\alpha,\alpha)e_\alpha.\eqno(4.25)
$$
If we restrict the relation (4.25) to the real simple roots, then it
determines a unique element in the Cartan subalgebra of \g, namely
$$
H = {\textstyle {1\over 2}}\sum_{\alpha > 0}[E_\alpha,E_{-\alpha}].
$$
Therefore, there is a unique
element in the extended Cartan subalgebra, of the form
$$
\hat H = H + gd,\eqno(4.26)
$$
such that (4.25) holds for the affine root $e_0$ as well. The integer
$g$ is the dual Coxeter number of \g.  The redefined operator is
$$\eqalign{
\zeta \Phi(z) &:= (\hat H \otimes 1)\Phi(z) + \normal{(\hat r^{\rm t} + \hat r)\Phi(z)}\cr &= (\hat H
\otimes 1)\Phi (z) +
\hat r^{\rm t}(z)\Phi(z) + (\hat r^i(z) \otimes 1)\Phi(z)\hat r_i.\cr}\eqno(4.27)
$$ 
Notice that we did not actually use (4.24); instead we defined $\hat
H$ as the element of $\hg$ that has the same commutator with
$e_\alpha$ as (4.23).  This makes it plausible that the term that was
dropped is a scalar and that covariance is preserved. Indeed we have
the
\b
\no{\bf Proposition 4.2.} {\it The operator product $\zeta\Phi(z)$ 
defined in (4.27) is covariant; that is, if $\Phi$ is an intertwiner
then so is $\zeta\Phi$.}
\b
\no {\it Proof.} The coproduct is that of the classical limit, 
$\Delta(x) = x \otimes 1 + 1 \otimes x$ for $x \in \hg$.
$$
\eqalign{
\Delta(x)\zeta \Phi(z) - \bigl(\zeta\Phi(z)\bigr) x &= 
([x,\hat H] \otimes 1)\Phi(z) 
+ [\Delta(x),\hat r^{\rm t}(z)]\Phi(z)\cr & + ([x,\hat r^i(z)] \otimes
1)\Phi(z)\hat r_i + (\hat r^i(z) \otimes 1)\Phi [x,\hat
r_i\,].\cr}\eqno(4.28)
$$
Suppose first that $x \in \g$; then in terms 2, 3, 4 only the zero
modes contribute. The sums over the degree are now finite and
$$
\eqalign{
\Delta(x)\zeta \Phi(z) - \zeta\Phi(z) x &= ([x, H] \otimes 1)\Phi(z) + 
([x, r^i(0)r_i(0)] \otimes 1)\Phi(z),\cr}\eqno(4.29)
$$
which vanishes in view of (4.26). We shall verify that (4.28) holds
for $x = e_0$.  Besides (4.29) there are additional terms that arise
from the extension term in the commutation relations, others that
arise from the fact that $e_0$ does not commute with the degree
operator, and finally the more subtle contributions that come from the
fact that the degree of $[e_0,y]$ is shifted by 1 from that of $y$:
\ve
$$
([e_0,\hat H] \otimes 1)\Phi = (-ge_0\otimes 1)\Phi + 
([e_0, H] \otimes 1)\Phi,
$$
$$
\eqalign{
[\Delta(e_0),\hat r^{\rm t}(z)]\Phi &= [\Delta(e_0), \hat r^{\rm
t}(0)]\Phi +
\sum_{n > 0}[\Delta(e_0),L_a^n \otimes L_a^{-n}]\Phi\cr 
&= [\Delta(e_0), \hat r^{\rm t}(0)]\Phi +  (L_a^1 \otimes [e_0,L_a^{-1}])  \Phi,   \cr} 
$$
$$
\eqalign{
([e_0,\hat r^i(z)]\otimes 1)\Phi\hat r_i &= 
([e_0, \hat r^i(0)]\otimes 1)\Phi \hat r_i + 
\sum_{n > 0}([e_0,L_a^{-n}]\otimes 1)\Phi L_a^n, \cr}
$$
$$
(\hat r^i(z) \otimes 1)\Phi [e_0,\hat r_i\,] =( \hat r^i(0) \otimes
1)\Phi [e_0, \hat r_i\,] +\sum_{n > 0}(L_a^{-n}\otimes
1)\Phi[e_0,L_a^n].
$$
The two infinite sums almost cancel, leaving only the first term of
the second one. The sum of the last two expressions is
$$
[\Delta(e_0),\hat r(0)]\Phi + ([e_0,r^i(0)r_i(0)] \otimes 1)\Phi +
([e_0,L_a^{-1}] \otimes 1)\Phi L_a^1.
$$
Adding the second expression we obtain
$$\eqalign{ [\Delta(e_0),\hat r(0) + \hat r^{\rm t}(0)]\Phi &+
 ([e_0,L_a^{-1}] \otimes L_a^1)\Phi + (L_a^1 \otimes
 [e_0,L_a^{-1}])\Phi\cr &+ ([e_0,r^i(0)r_i(0)] \otimes 1)\Phi +
 ([e_0,L_a^{-1}]L_a^1 \otimes 1)\Phi.
\cr}
$$
The first three terms cancel exactly and we have 
$$
\eqalign{
\Delta(e_0)\zeta \Phi(z) - \zeta\Phi(z)e_0 &= 
-g(e_0\otimes 1)\Phi + ([e_0, H] \otimes 1)\Phi\cr
& + ([e_0,r^i(0)r_i(0)] \otimes 1)\Phi + ([e_0,L_a^{-1}]L_a^1 \otimes
1)\Phi.
\cr}
$$
Terms two and three cancel as in (4.29) and the proposition is proved
when we verify that, in the evaluation module, $ [e_0,L^{-1}_a]L^1_a =
[e_0,L_a]L_a = ge_0, $ and repeat the calculation with $e_0$ replaced
by $e_{-0}$.
\b
\no {\it Normalization.} 

Returning  to (4.27) we put the degree operator into evidence:
$$
\zeta\Phi(z) = (k+g)z{{\rm d}\over{\rm d}z}\Phi + (H \otimes 1)\Phi(z)  +   
r^{\rm t}(z)\Phi(z) + \bigl(r^i(z) \otimes
1\bigr)\Phi(z)r_i.\eqno(4.30)
$$
Again we fix the grading of the intertwiner as in (4.21), but now with
$(\mu|\nu,\pi)$ replaced by
$$ 
(\mu|\nu) = {A(w)\over k+g},\quad  
A := \varphi(.,v_0)  + \varphi(v'_0,.) + H  = 
{\textstyle {1\over 2}}\bigl(C(\mu) - C(\nu)\bigr), 
\eqno(4.31)
$$
\vskip1mm\no so that the operator form of the Knizhnik-Zamolodchikov 
equation takes the form
$$
\zeta\Phi(z) = 0.\eqno(4.32)
$$
Note that this makes the grading of $\Phi$ independent of the choice
of evaluation module; this grading/normalization is thus ``universal".
\b
\no {\it Remark 4.3.} In view of the interpretation of the quantum 
field $\Phi(z)$ as an intertwiner for highest weight affine Kac-Moody
modules, the appearance of the {\it rational} r-matrix in the original
KZ equation (4.22) has always seemed somewhat mysterious. The mystery
is deepened by the discovery [K] that the monodromy associated with
the solutions yields a representation of $U_q($\g).  The alternative,
to use the polarization based on the decompostion $\zeta = \hat r +
\hat r^{\rm t}$, was first suggested in [FR]; it seems to be more
natural. However, $\Phi$ is defined as an intertwiner of Kac-Moody
modules, with the classical coproduct; it knows nothing about
r-matrices.  Normal ordering is an example of additive
renormalization, or ``subtraction", necessary only if the ordinary
product is ill defined.  Any two polarizations that eliminate the
divergent term by subtracting a scalar (that is; without compromising
covariance) are equivalent, and one is not more natural than the
other, in the present context at least. The fact that renormalization
is required is revealed by the fact that the subtracted term, the last
term in (4.15), is divergent. It is related to the fact that the
classical r-matrix has a pole at $z_1/z_2 = 1$.

\b
\no {\it Remark 4.4.} The appearance of the factor $k + g$ as a 
coefficient of the degree operator in both versions is justified by
covariance, as is indeed implied by the proof of Proposition 4.1 in
[FR]. The term $ (H \otimes 1)\Phi(z)$ has exactly the same
origin. Perhaps it should be pointed out that the concept of
``covariance" that is evoked in this Section is quite distinct from
the covariance under twisting that is alluded to in the title of the
paper and in from Section 6 forward.

\bb

\no {\bf {\steptwo 5. The quantum KZ equation.}}
\b\no
Here we shall make sense of Eq.(4.2),
 $$
(Z - \langle Z \rangle)\Phi(z) = 0,\quad Z = RR^{\rm t},
$$ 
in the quantized Kac-Moody algebra, to recover the q-KZ equation of
Frenkel and Reshetikhin [FR].
 
The action of $R^{\rm t}\Phi$ on $V_{\mu,k}$ is well defined (in terms
 of formal series), since both $\Phi$ and $R^{\rm t}$ act by
 degree-decreasing operators in the second space (Proposition 3.1 and
 Eq.(3.9)), but the subsequent action of $R$ is not. We therefore
 investigate the effect of normal ordering. Thus if
$$
R = R^i \otimes R_i,
$$
we set (tentatively) 
$$
\normal{Z\Phi} = (R^i \otimes 1)\Psi R_i,\quad \Psi:= R^{\rm t}\Phi\quad (?), 
$$
and try to prove that the operator $Z: \Phi \mapsto \normal{Z\Phi}$ is
invariant; that is, that it commutes with the coproduct. In view of
the intertwining property of $\Phi$ this is the same as
$$
\Delta(x)\normal{Z\Phi} =\normal{Z\Phi}x\quad (?).
$$
Attempts to verify this equation leads to
\b
\no {\bf Proposition 5.1.} {\it Let $\hat H$ be the element in 
the Cartan subalgebra ${\cal A}'_0$ of ${\cal A}'$ with the property
$$
q^{\hat H}e_\alpha q^{-\hat H} =
q^{\varphi(\alpha,\alpha)}e_\alpha,\eqno(5.1)
$$
and define the normal-ordered product $\normal{Z\Phi}$ by
$$
\normal{Z\Phi}:= (\hat R^i \otimes 1)\Psi R_i,\quad \Psi:= R^{\rm t}\Phi,
\quad \hat R^i :=  R^i q^{\hat H}.\eqno(5.2) 
$$
Then
$$
\Delta(x)\normal{Z\Phi} = \normal{Z\Phi}x,\quad \forall x \in 
{\cal A}'.\eqno(5.3)
$$
}
\b
\no {\it Proof.} We begin with
$
R\,\Delta'(e_\alpha) =  \Delta(e_\alpha)R; 
$
that is
$$
 (R^i \otimes R_i)(e_\alpha \otimes1 + q^{\varphi(\alpha,.)} 
\otimes e_\alpha) = (1 \otimes e_\alpha + e_\alpha \otimes
q^{\varphi(\alpha,.)})(R^i \otimes R_i),
$$
and thus
$$
R^i \otimes R_ie_\alpha = -R^ie_\alpha q^{-\varphi(\alpha,.)} \otimes
R_i + R^iq^{-\varphi(\alpha,.)} \otimes e_\alpha R_i + e_\alpha
R^iq^{-\varphi(\alpha,.)} \otimes q^{\varphi(\alpha,.)}R_i\,,
$$
which gives us
$$\eqalign{
\normal{Z\Phi}e_\alpha = (R^iq^{\hat H} &\otimes 1)\Psi R_ie_\alpha = 
-(R^ie_\alpha q^{-\varphi(\alpha,.)}q^{\hat H}\otimes 1)\Psi R_i\cr &+
(R^iq^{-\varphi(\alpha,.)}q^{\hat H} \otimes 1)\Psi e_\alpha R_i +
(e_\alpha R^iq^{-\varphi(\alpha,.)}q^{\hat H} \otimes 1)\Psi
q^{\varphi(\alpha,.)}R_i\,.
\cr}
$$
Using the intertwining property of $\Psi$ we convert the last two
terms to
$$
(R^i q^{\hat H} \otimes e_\alpha)\Psi R_i +
(R^iq^{-\varphi(\alpha,.)}q^{\hat H}e_\alpha \otimes 1)\Psi R_i +
(e_\alpha R^i q^{\hat H}\otimes q^{\varphi(\alpha,.)})\Psi R_i.
$$
As for the first term, we shift the operator $e_\alpha$ to the right;
since $e_\alpha$ commutes with $\hat H -
\varphi(\alpha,.)$   we get the required cancellation and the   
result is
$$
\normal{Z\Phi}e_\alpha =  
 ( \hat R^i \otimes e_\alpha)\Psi R_i + (e_\alpha\hat R^i \otimes
q^{\varphi(\alpha,.)})\Psi R_i = \Delta(e_\alpha)(\hat R^i \otimes
1)\Psi R_i.
$$
 In the classical limit the $q$-factor in (5.2) produces the $\hat
H$-term in Eq.(4.27).  A similar calculation with $e_{-\alpha}$
completes the proof of Proposition 5.1, and we have an independent
confirmation of the covariance of (4.27).
\b
\no {\it Normalization.}

Our next task is to pull out the degree operator. Since the first
space is an evaluation module, on which the central element $c$ is
zero, the degree operator appears only in the first factors of $R$ and
$R^{\rm t}$, as $z{{\rm d}\over{\rm d}z}$.  We define $L^{\mp}$ by
$$
\hat  R^i(z)\otimes R_i  =  q^{(1-u)kd + gd}\bigl({\rm Ad}(q^{-gd})\otimes 
1\bigr)L^-(z),\quad R^{\rm t}(z) = q^{ukd}L^+(z),\eqno(5.4)
$$
where $d = z{{\rm d}\over{\rm d}z}$ acts in the evaluation module and
Ad$(x)y = xyx^{-1}$.  Objects denoted by the letter $L$ (with
ornamentation) do not contain $d$. We also need the expansions
$$
L^-(z) = L^{-i}(z) \otimes L^-_i,
\quad L^+(z)=L^{+i}(z)\otimes L^+_i.
$$
Now
$$
\eqalign{
\normal{Z\Phi(z)} & = q^{(1-u)kd + gd}\bigl({\rm Ad}(q^{-gd})L^{-i}(z) 
\otimes 1\bigr)(q^{ukd} \otimes
1)L^+(z)\Phi(z)L_i^-\cr & = q^{(k +g)d}\bigl(L^{-i}(q^{-g-uk}z)
\otimes 1\bigr)L^+(z)\Phi(z)L^-_i\cr &=:
q^{(k+g)d}\normal{L(z)\Phi(z)}.
\cr}\eqno(5.5)
$$
Thus, the q-KZ equation for $\Phi$:
$$
 \Phi(q^{-k-g}z) = \normal{L(z)\Phi(z)} = 
\bigl(L^{-i}(q^{-g-uk}z) \otimes 
1\bigr)L^+(z)\Phi(z)L^-_i.\eqno(5.6)
$$
Here we have fixed $\langle Z \rangle = 1$. It means that the grading
of $\Phi$ is so chosen that, on any weight vector $w \in V$ and the
highest weight vectors $v_0 \in V_{\mu,k}, v'_0 \in V_{\nu,k}$, we
have
$$
  {(w\otimes v_0',\Phi(q^{-k-g}z)v_0)\over (w\otimes v_0',\Phi(z)v_0)}
= q^{(\mu|\nu)},
$$
with $(\mu|\nu)$ as in (4.31).

The other intertwiner, $\Psi \propto R^{\rm t}\Phi$, satisfies $\Psi x
= \Delta'(x)\Psi$ and $\normal{Z'\Psi} = \Psi$ where $Z' := R^{\rm
t}R$. We find
$$
\normal{Z'\Psi(z)}= R^{\rm t}(\hat R^i \otimes 1)\Psi(z)R_i = 
q^{ukd}L^+(z)q^{(k-uk+g)d}\bigl(L^{-i}(q^{-g}z) \otimes
1\bigr)\Psi(z)L_i^-,\eqno(5.7)
$$
and thus,
 the q-KZ equation for $\Psi$:
$$
\Psi(q^{-k-g}z) =\normal{L'(z)\Psi(z)}= 
L^+(q^{-g-k+uk}z)\bigl(L^{-i}(q^{-g}z) \otimes
1\bigr)\Psi(z)L_i^-.\eqno(5.8)
$$

\bb
\no{\bf {\steptwo  6. Hopf twisting.}}
\b\no
It is remarkable that the elliptic quantum group can be viewed as
deformation of the trigonometric quantum group.  The deformation does
not affect the algebraic structure, which remains that of a quantized,
affine Kac-Moody algebra. Only the coproduct distinguishes the
elliptic case from the trigonometric one. The deformation is
implemented by a twist in the category of Hopf algebras (this section)
or quasi Hopf algebras (next section). The full elliptic quantum group
is quasi Hopf; it becomes Hopf on the quotient by the ideal generated
by the center.  In this Section we investigate the effect of twisting
on the intertwiners and on the KZ equation, in the quantum case where
the relationship between the intertwiner and the R-matrix is more
clear.
\b
\no {\bf Definition 6.1.} A formal Hopf deformation of    
  a standard R-matrix $R$ is a formal power series
$$
R_\epsilon = R + \epsilon R_1 + ...~,
$$
that satisfies the  Yang-Baxter relation to each order in $\epsilon$.  
\b
It turns out [Fr1] that the deformations of greatest interest have the
form of a twist.
\b
\no {\bf Theorem 6.2.} {\it Let $R$ be the R-matrix, $ \Delta$  
the coproduct, of a coboundary Hopf algebra ${\cal A}',$ and $F \in
{\cal A}'\otimes {\cal A}'$, invertible, such that
$$
\bigl((1 \otimes \Delta_{21}) F\bigr)F_{12} = \bigl((\Delta_{13} \otimes 1)
F\bigr) F_{31}.
\eqno(6.1)
$$
 Then
$$
\tilde R := (F^{\rm t})^{-1} R F\eqno
$$
(a) satisfies the Yang-Baxter relation and (b) defines a Hopf algebra
$\tilde {\cal A}$ with the same product and with coproduct
$$
\tilde \Delta = (F^{\rm t})^{-1} \Delta F^{\rm t}.\eqno$$
}
\b
\no This is a result of Drinfel'd [D2]; a detailed proof was given in [Fr1].
\b
We say that a deformation $R_\epsilon$ of a standard R-matrix $R$ is
implemented by a {\it twistor} $F_\epsilon$ if there is a formal power
series
$$
F_\epsilon = 1 + \epsilon F_1 + ...
$$
that satisfies (6.1) to each order in $\epsilon$ and
$$
R_\epsilon = (F_\epsilon^{\rm t})^{-1}RF_\epsilon.\eqno(6.2)
$$
In this case the deformed R-matrix intertwines a deformed coproduct,
$$
R_\epsilon\Delta_\epsilon' =
\Delta_\epsilon R_\epsilon,~~~ \Delta'_\epsilon := 
F_\epsilon^{-1}\Delta'F_\epsilon.\eqno(6.3)
$$
 
 Known solutions of (6.1) have the following structure [Fr1]. We need
a pair of subalgebras $\Gamma_1,\Gamma_2$ of ${\cal A}'=\hg$,
generated by sets $\hat \Gamma_i \subset \{e_\alpha\}_{\alpha \in N}$,
and a diagram isomorphism $\tau: \hat\Gamma_1 \rightarrow
\hat\Gamma_2$. A deformation exists when the parameters of ${\cal A}'$
satisfy the following condition,
$$
\varphi(\sigma,.) + \varphi(.,\tau\sigma) = 0,~~~\sigma \in \hat\Gamma_1.
$$
Note that $e_{\tau\sigma}$ is defined only if $e_\sigma \in \hat
\Gamma_1$.  Then there is a cocycle $F_\epsilon$ of the form
$$
\eqalign{
F_\epsilon = \prod_{m \geq
1}F_\epsilon^m:=F_\epsilon^1F_\epsilon^2...F_\epsilon^m...,&\quad
F_\epsilon^m =
\sum_{(\sigma)}
\epsilon^{mn}F_{(\sigma)}^{m(\rho)} f_{\sigma_1}...f_{\sigma_n} \otimes
f_{-\rho_1}...f_{-\rho_n},\cr f_\sigma :=
q^{-\varphi(\sigma,.)}e_\sigma, &\quad f_{-\rho} :=
e_{-\rho}\,q^{\varphi(.,\rho)},}
\eqno(6.4)
$$
where the sum is over all $(\sigma) = \sigma_1,...,\sigma_n$, and all
permutations $(\sigma')$ of $(\sigma)$, such that $\rho_i =
\tau^m\sigma'_i$ is defined. We take $F_{(\sigma)}^{m(\rho)} = 1$ when
the set $(\sigma)$ is empty.

Note that the family of deformation of this type is large enough to
contain the quantization of all the classical Lie bialgebras
classified by Belavin and Drinfel'd, with r-matrices of constant,
trigonometric and elliptic type.  Two cases need to be distinguished.

(a) {\it Finite twisting} is by definition the case when there is $k$
such that for all $\sigma$, $\tau^k\sigma \notin \hat\Gamma_1$; then
$\hat\Gamma_1,\hat\Gamma_2$ are distinct and the product over $m$ is
finite.

(b) {\it Elliptic twisting.} The only other possibility (see [Fr1],
 Section 16) is that ${\cal A}' = \widehat {\sl(N)}$ and $\Gamma_1 =
 \Gamma_2$ is generated by all the simple roots. This section deals
 with twisting in the category of Hopf algebras; elliptic twisting
 within the context of Hopf algebras implies [Fr2] that we drop the
 central extension and descend to loop algebras. The full elliptic
 Kac-Moody algebra is quasi Hopf and will be discussed in the next
 section.

The deformed R-matrix and coproduct are
$$
R_\epsilon = (F_\epsilon^{\rm t})^{-1}RF_\epsilon,\quad
\Delta_\epsilon = (F_\epsilon^{\rm t})^{-1}\Delta F_\epsilon^{\rm
t}.\eqno(6.5)
$$
Here are some products that seem ill defined; thus $R$ has
degree-increasing operators in the second space, where $F_\epsilon$
has degree-decreasing operators. This problem can be handled in a
general way by adopting an interpretation that is quite natural in
deformation theory. One notes that $F_\epsilon$ is a formal power
series in the deformation parameter $\epsilon$. One interprets all the
operators this way; then the problem reduces to making sure that the
coefficients are well defined. Indeed, to any fixed order in
$\epsilon$, the product $RF_\epsilon$ is, in the second space, a power
series in the operators $e_\alpha$ multiplied by a polynomial in the
other generators.

It is, nevertheless, of some interest to determine whether
singularities arise as one assigns a value to $\epsilon$ and attempts
to sum up the deformation series. In this respect cases (a) and (b)
are quite different.

(a) Finite twisting. The sum in (6.4) becomes finite when projected on
a finite dimensional representation in either one of the two
spaces. Infinite sums will appear if both representations are
infinite, but there is a finite number of terms with fixed weight;
therefore no infinite, purely numerical series will appear. Infinite
sums with operator coefficients are beyond (our power of) analysis in
the general case, and of no immediate concern to us. The value of
$\epsilon$ is basis dependent; the only distinct possibilities are
$\epsilon = 0,1$.

(b) We note that the range of $\epsilon$ is in this case $|\epsilon|
<1$.  Here the situation is more delicate, and of some interest.
Under twisting, the Casimir element $Z$ suffers an equivalence
transformation
$$
Z_\epsilon = (F_\epsilon^{\rm t})^{-1}ZF_\epsilon^{\rm t},\eqno(6.6)
$$
and one expects that an intertwiner $\Phi_\epsilon$, satisfying
$$
\Delta_\epsilon(x)\Phi_\epsilon = \Phi_\epsilon x,\quad x 
\in {\cal A}',\eqno(6.7)
$$
may be expressed as $\Phi_\epsilon = (F_\epsilon^{\rm
t})^{-1}\Phi$. However, $F_\epsilon^{\rm t}$ has a structure similar
to that of $R$, with degree-increasing operators in the second space,
and we must consider the possibility that normal ordering may be
required. In fact probably not, but since we have not proved this, we
shall switch our attention to the other intertwiner.

We consider instead the alternative intertwiner $\Psi$, and the
alternative Casimir operator $Z'$ that commutes with $\Delta'(x)$,
namely
$$
Z' = R^{\rm t}R,\quad Z'\Delta'(x) = \Delta'(x)Z',\quad (Z' - 1)\Psi = 
0.\eqno(6.8)
$$
We have
$$
Z'_\epsilon = F_\epsilon^{-1}Z'F_\epsilon, \quad \Psi_\epsilon =
F_\epsilon^{-1}\Psi.\eqno(6.9)
$$
The operator product $F_\epsilon^{-1}\Psi$ is well defined as an
operator on $V_{\mu,k}$.  The intertwining property of
$\Psi_\epsilon$, namely
$$
\Delta'_\epsilon(x)\Psi_\epsilon = \Psi_\epsilon x,\eqno(6.10)
$$
is therefore in the clear.

We define the operator product
$Z'_\epsilon\Psi_\epsilon(z)$. Formally,
$$
Z'_\epsilon\Psi_\epsilon = R^{\rm t}_\epsilon R_\epsilon\Psi_\epsilon
= F_\epsilon^{-1}Z'\Psi,
$$
and this too is well defined, provided we define the untwisted product
as in (5.7); that is
$$
Z'\Psi \rightarrow \normal{Z'\Psi} = R^{\rm t}( R^iq^{\hat H} \otimes
1)\Psi R_i.
$$
The equation satisfied by the twisted correlation function is
$(Z'_\epsilon-1)\Psi_\epsilon = 0$ or more precisely
\b
\no {\bf Definition 6.3.} The twisted q-KZ equation is the following 
equation for the twisted intertwiner operator,
$$
\Psi_\epsilon(q^{-k-g}z) =  F_\epsilon^{-1}(q^{-k -g}z) L^+(zq^{-g-k+uk})(L^{-i}(zq^{-g}) \otimes
1\bigr)F_\epsilon(z)\Psi_\epsilon(z)L^-_i.\eqno(6.11)
$$
 
 It should be noted that the polarization used is the same as before
deformation.  To justify this we repeat that the definition of the
intertwining operators is independent of normal ordering conventions,
normal ordering is relevant only when the ordinary product does not
exist, it is required to be well defined and covariant, nothing
more. Of course, it is also true that, if $\Psi_\epsilon$ is defined
as in (6.9), then (6.11) is equivalent to (5.8).

The top matrix element  of $\Psi_\epsilon$ is
$$
(v'_0,\Psi_\epsilon v_0) = (v'_0,F_\epsilon^{-1}\Psi  v_0) = 
(v'_0,\Psi  v_0).\eqno(6.12)
$$
This shows that, in a complete description of, say, the eight-vertex
model, both periodic and non-periodic functions appear. We had naively
expected to encounter nothing but elliptic functions, that ``the
eight-vertex model lives on the torus".

Having thus discarded a prejudice, we are comfortable with the
continued use, in the twisted case, of the original polarization based
on the standard trigonometric R-matrix. The alternative of defining a
normal-ordered product such that
$$
R_\epsilon\Phi_\epsilon = (R_\epsilon^i \otimes 1)\Phi_\epsilon R_{\epsilon i}
$$
is entirely redundant.
 
Another idea is to replace matrix elements by traces, as suggested by
Bernard [Ber] and in [FR]. However, since we know that the elliptic
quantum group, as an algebra, is the same as the standard quantum
group (that is, a Kac-Moody algebra), there seems to be no reason to
take less interest in the highest weight matrix elements in the
elliptic case. Continuity of physics also suggests that we continue to
work with the usual module structure, as was argued in [JMN], Section
4. Trace functionals are interesting in themselves, but there seems to
be no reason to neglect the matrix elements.

  The intertwiners of Kac-Moody modules, and the solutions of the KZ
equation, know nothing about r-matrices. For all that we may derive
different versions of the equation, the solutions remain the same. To
base the polarization on the R-matrix is not an imperative; more
important is to adopt a workable definition that gives a meaning to
the objects of interest; to wit, matrix elements of intertwiners.

In the setting of conformal field theory twisting does not affect the
quantization paradigm, but it does change the quantum fields (the
intertwiners) and their operator product expansions.

We shall need to know the twistor $F_\epsilon$. It is determined,
uniquely, by the recursion relations [Fr1]
$$
[1 \otimes f_\rho,F_\epsilon^m] = 
\epsilon^m\bigl( F_\epsilon^m(f_{\tau^{-m}\rho} \otimes
q^{-\varphi(\rho,.)}) -
 (f_{\tau^{-m}\rho} \otimes q^{\varphi(.,\rho)})F_\epsilon^m\bigr),\eqno(6.13)
$$ 
with the initial conditions
$$
F_\epsilon^m = 1 - \epsilon^m\sum_{(\rho)}f_{\tau^{-m}\rho} \otimes
f_{-\rho} + ...\,.\eqno(6.14)
$$
These equation were solved in a special case, and used to calculate
the elliptic R-matrix of $\widehat{\sl(2)}$ in the fundamental
representation [Fr1].  Later, we shall exploit the similarity between
this relation and the recursion relation (2.4) for the universal
R-matrix.

\bb
\no {\bf {\steptwo 7. Quasi Hopf twisting.}}
\b\no
We are interested in the elliptic quantum groups, in the sense of
Baxter [Ba] and Belavin [Be].  This takes us out of the framework of
Hopf algebras, but just barely so. The special nature of these quasi
Hopf algebras is that they become Hopf algebras at level zero; that
is, on the quotient by the ideal generated by the center.

Quasi Hopf deformations are constructed in the same way as Hopf
deformations, except that the element $F_\epsilon$ need not satisfy
the cocycle condition (6.1). The deformed R-matrix and coproduct are
given by (6.5), but the former no longer satifies the Yang-Baxter
relation and the latter is not coassociative, in general.

 If $F^{m(\rho)}_{(\sigma)}$ are the coefficients of the elliptic Hopf
 twistor in (6.4), then the elliptic quasi Hopf twistor has the form
 [Fr2]
$$
  F_\epsilon = \prod_{m = 1,2,...}  F_\epsilon^m,\quad F_\epsilon^m =
\sum_{(\sigma)}\epsilon^{nm}F^{m(\rho)}_{(\sigma)}
f_{\sigma_1}...f_{\sigma_n} \otimes f_{-\rho_1}...f_{-\rho_n}Q(m,\rho),\eqno(7.1) 
$$
where $Q(m,\rho) \in {\cal A}'_0 \otimes {\cal A}'_0$ and ${\cal
A}'_0$ is the Cartan subalgebra of the quantized Kac-Moody algebra
${\cal A}'$. This factor is equal to unity in the Hopf case, and (7.1)
then reduces to (6.4). The $F_\epsilon$-twisted algebra is a Hopf
algebra when the parameters satisfy the condition
$$
\varphi(\sigma,.) + \varphi(.,\tau\sigma) = 0,~~~\sigma \in \hat\Gamma_1,
$$
where now $\tau$ is the cyclic diagram automorphism that takes each
simple root of $\sl(N)$ to its neighbour. This condition can be
satisfied on the loop algebra (when $c \mapsto 0$). We are interested
in the full Kac-Moody algebra ($c \neq 0$); in that case the best that
can be done is to choose parameters such that
$$
\varphi(\sigma,.) + \varphi(.,\tau\sigma) = [(1-u)\delta_\sigma^0 + u\delta_{\tau\sigma}^0]c.
\eqno(7.2)
$$ 
This algebra is what we mean by ``elliptic quantum group in the sense
of Baxter and Belavin"; it is a quasi Hopf algebra of a particularly
benevolent type, where the deviation from coassociativity is confined
to the center.

Instead of the cocycle condition (6.1) we now have 
$$
\bigl(({\rm id} \otimes \Delta_{21})F_\epsilon\bigr)F_{\epsilon 12} = 
\bigl((\Delta_{13}\otimes {\rm id})F_\epsilon\bigr)F_{\epsilon 31,2},\eqno(7.3)
$$
where $F_{\epsilon ij,k}$ is an extension of $F_{\epsilon ij}$,
supported on the center, to the $k$'th space. In the case of interest,
when we are dealing with modules with fixed level $c
\mapsto k$, this amounts to a  modification of the coefficients in $F_{\epsilon ij}$. From (7.3)
one gets the Cartan factors
$Q(m,\rho)$ [Fr2] and the recursion relation
$$
[1 \otimes f_\rho,F_\epsilon^m] = \epsilon^m\bigl(
F_\epsilon^m(f_{\tau^{-m}\rho} \otimes q^{-\varphi(\rho,.)}) -
(f_{\tau^{-m}\rho} \otimes
q^{\varphi(.,\rho)})F_\epsilon^m\bigr)Q(m,\rho),\eqno(7.4)
$$ 
with the initial conditions
$$
F_\epsilon^m = 1 - \epsilon^m\sum_\rho(f_{\tau^{-m}\rho} \otimes
f_{-\rho})Q(m,\rho) + ...\,.\eqno(7.5)
$$
The solutions will be given later.
Once $F_{\epsilon 12}^{m}$ is known, $F^m_{\epsilon 12,3}$ is obtained by means of the
substitution 
$$
1 \otimes c \mapsto 1 \otimes \Delta(c).\eqno(7.6)   
$$

\ve
\no{\bf {\steptwo 8. Correlation Functions.}}
\b\no
The main objects of interest, in conformal field theory as well as in
the study of statistical models, are the correlation functions. In
their simplest form they are matrix elements of products of
intertwiners,
$$\eqalign{
&f_{v'v}(z_1,...,z_N) = \langle v',\Phi(z_1)...\Phi(z_N)\,v\rangle \in V_1(z_1) \otimes ...\otimes V_N(z_N),\cr &
g_{v'v}(z_1,...,z_N) = \langle v',\Psi(z_1)...\Psi(z_N)\,v
\rangle \in V_1(z_1) \otimes ...\otimes V_N(z_N).\cr}
\eqno(8.1)
$$
 Here $\Phi(z_p)$ and $\Psi(z_p)$ are intertwiners between highest
weight modules,
$$
\Phi(z_p), \Psi(z_p): V_{\mu_p,k} \rightarrow V_p(z_p) 
\otimes V_{\mu_{p-1},k},~~~ p = 1,...,N,
$$
with $\{V_p(z_p)\}$ a set of evaluation modules, and $v \in
V_{\mu_N,k}$, $v'
\in V_{\mu_0,k}$. These ``functions" are formal, $V_1 \otimes 
... \otimes V_N$-valued 
series in $N$ distinct variables.
\b
\no {\it Classical Correlation Functions.}

We begin with the classical case and the polarization (4.18), 
$$
J_a = J_a^+ + J_a^- = \sum_{n > 0}z^n L_a^{-n} + \sum_{n \geq
0}z^{-n}L_a^n,
$$
and the normalization that leads to (4.22):
$$
 (k+g)z{{\rm d}\over{\rm d}z}\Phi(z) + L_a \otimes J_a^+\Phi(z) + (L_a
\otimes 1)\Phi(z)J^-_a = 0.
$$
Then for any $p \in \{1,...,N\}$,
$$
\eqalign{
&(k+g)z_p{{\rm d}\over{\rm d}z_p}f_{v'v}(z_1,...,z_N) =\cr
&-L_a^{(p)}\langle v',...\Phi(z_{p-1}) J_a^+(z_p)\Phi(z_p)... v\rangle
- L_a^{(p)}\langle v',...\Phi(z_{p-1})
\Phi(z_p)J_a^-(z_p)... v\rangle.
\cr}\eqno(8.2)
$$
Here $L_a^{(p)}$ denotes the action of $L_a$ in $V_p$. Suppose now
that the vectors $v_0$ and $v_0'$ are highest weight vectors of the
respective highest weight modules. The intertwiners satisfy
$[L^n_a,\Phi(z_p)] = -L_a^{(p)}z_p^n\Phi(z_p)$; this allows us to
permute $J^+$ through to the left, where it dies on the highest weight
vector, and to permute $J^-$ towards the right, where only the zero
modes survive, to contribute the last term in
$$
\eqalign{
&(k+g)z_p{{\rm d}\over{\rm d}z_p}f_{v_0'v_0}(z_1,...,z_N) = 
-\sum_{1\leq q<p}~\sum_{n > 0}
({z_p\over z_q})^nL_a^{(q)}L_a^{(p)}\langle v_0',...\Phi(z_{p-1})
\Phi(z_p)... v_0\rangle\cr &  + 
\sum_{N\geq q > p}~\sum_{n\geq 0}({z_q\over z_p} )^n
L_a^{(p)}L_a^{(q)}\langle v'_0,...\Phi(z_{p-1})
\Phi(z_p) ... v_0\rangle + L_a^{(p)}\langle v',...\Phi(z_{p-1})\Phi(z_p) 
... L_a  v_0\rangle.\cr}
$$
Hence
$$
(k+g){{\rm d}\over{\rm d}z_p}f_{v_0'v_0}(z_1,...,z_N) = 
\sum_{q \neq p}{1 \over z_p - z_q}L_a^{(p)}L_a^{(q)}f_{v_0'v_0}(z_1,...z_N),
\quad p = 1,...,N,\eqno(8.3)
$$
where $q$ takes the values $1,...,N+1$, $z_{N+1} = 0$, and
$L_a^{(N+1)}$ acts on $v_0$.  The last expression must be supplemented
by the instruction
$$
{1\over z_p - z_q} := \cases{ (1/ z_p)\displaystyle\sum_{n\geq 0}
(z_q/z_p)^n, & $ q > p$, \cr (-1/z_q)\displaystyle\sum_{n\geq
0}(z_p/z_q)^n, & $q < p$.\cr}\eqno(8.4)
$$
The domain of convergence is thus $|z_1| > |z_2| > ...> |z_{N+1}| =
0$.

In the simplest, nontrivial case $N = 1$.  Projecting on a vector $w
\in V$ we get
$$
 (k+g){{\rm d}\over{\rm d}z_1}f^w_{v_0'v_0}(z_1) = 
{c \over z}f^w_{v_0'v_0}(z_1), 
\quad c = {\langle w \otimes v',{\cal C}\Phi\, v \rangle\over \langle w 
\otimes v', \Phi v\rangle} = {\textstyle {1\over 2}}\bigl(C(\nu) - C(\mu) -
C(\pi)\bigr),
$$
which simply reflects the choice of grading of $\Phi$.  The case $N =
2$ is not much more complicated; the equations are
$$
 (k+g){{\rm d}f\over{\rm d}z_1} = {c_{12}f\over z_1 - z_2} +
{c_{13}f\over z_1},\quad (k+g){{\rm d}f\over{\rm d}z_2} =
{c_{12}f\over z_2 - z_1} + {c_{23}f\over z_2},\eqno(8.5)
$$
where $c_{ij} = L_a^{(i)}L_a^{(j)}$ and ``3" refers to the source
space.  In the case of $\widehat {\sl(2)}$ and fundamental evaluation
modules it is a simple matter to work out the hypergeometric
solutions. The general structure of the solution was exploited by
Khono [K] and Drinfel'd [D2] to construct representations of the braid
group and examples of quasi Hopf algebras.

If instead we use the polarization of (4.27) we obtain from (4.30) and (4.31), 
on the vectors of highest weight,
$$\eqalign{&
(k + g) z_p{{\rm d}\over{\rm d}z_p}f_{v_0'v_0}(z_1,...,z_N) + 
\bigl(A_p+\sum_{q=1}^{p-1}r_{qp}-
\sum_{q=p+1}^Nr_{pq}\bigr)f_{v_0'v_0}(z_1,...,z_N) = 0,\cr
& A_p := (H + \varphi(v_0',.) + \varphi(.,v_0))_p,}
\eqno(8.6)
$$
for $p = 1,...,N$. When $N=2$,
$$
(k + g) z_1{{\rm d}\over{\rm d}z_1}f_{v_0'v_0}(z_1,z_2) + A_1
f_{v_0'v_0}(z_1,z_2) - r_{12}f_{v_0'v_0}(z_1,z_2) = 0,
$$
$$
(k + g) z_2{{\rm d}\over{\rm d}z_2}f_{v_0'v_0}(z_1,z_2) + A_2
f_{v_0'v_0}(z_1,z_2) + r_{12}f_{v_0'v_0}(z_1,z_2) = 0.
$$
The solutions are, of course, the same, up to normalization.
\b
\no {\it q-Deformed Correlation Functions.}

We turn to the q-KZ equation (5.6), 
$$
 \Phi(q^{-k-g}z) = \normal{L(z)\Phi(z)}. 
$$
For functions of the type (8.1) the implication is
$$\eqalign{T_pf_{v_0'v_0}(z_1,...,z_N) &:=
f_{v_0'v_0}(...,z_{p-1},q^{-k-g}z_p,z_{p+1},...) \cr & ~=
\langle v'_0, ...\Phi(z_{p-1})L^{-i}(z'_p)L^+(z_p)\Phi(z_p) 
L^-_{i}\Phi(z_{p+1})...
 v_0\rangle ,
\cr}$$
with $z' = q^{-g-uk}z$. More transparently,
$$T_pf_{v_0'v_0}(z_1,...,z_N) = 
\bigl[L^{-i}(z'_p)L^{+j}(z_p)\bigr] 
\langle v_0',...\Phi(z_{p-1})L^+_{j}\Phi(z_p) L^-_{i}\Phi(z_{p+1})...
 v_0\rangle. \eqno(8.7)
$$
  For $N = 2$,
$$\eqalign{ & f_{v_0'v_0}(q^{-k-g}z_1,z_2) = \bigl[L^{-i
}(z'_1)q^{-\varphi(v'_0,.)}\bigr]_1
\langle v_0', \Phi(z_1) L^-_i\Phi(z_2)  v_0\rangle ,\cr 
& f_{v_0'v_0}(z_1,q^{-k-g}z_2) = \bigl[q^{\varphi(.,v_0)+H}L^{+i
}(z_2) \bigr]_2
\langle v'_0, \Phi(z_1) L^+_i\Phi(z_2)  v_0\rangle.\cr}\eqno(8.8)
$$
We reduce this using the quasi triangularity conditions in the
Appendix. The final result is
$$\eqalign{ & T_1f_{v_0'v_0}(z_1,z_2) = R^{-1}_{12}({z_2\over
z_1}q^{g+k})q^{A_1}f_{v_0'v_0}(z_1,z_2) ,\cr &T_2f_{v_0'v_0}(z_1,z_2)
= q^{A_2} R_{12}({z_2\over z_1})f_{v_0'v_0}(z_1,z_2),\cr & A_i :=
\bigl(\varphi(v'_0,.) + \varphi(.,v_0)+H\bigr)_i,\quad i =
1,...,N.\cr}
\eqno(8.9)
$$
These two equations can be combined in two ways. The result is the
same in either case, in consequence of the fact that the operator $
A_1 +A_2 $ (the subscripts refer to the two evaluation modules)
commutes with $R_{12}$.  From the fact that the correlation function
is invariant for the action of the Cartan subalgebra in the four
spaces it follows in fact that we can replace
$$
A_1+A_2 \rightarrow {\textstyle{1\over 2}}\bigl(C(\mu) - C(\nu)\bigr).
$$
The result is that
$$
T_1T_2f_{v_0'v_0}(z_1,z_2) = q^{A_1+A_2}f_{v_0'v_0}(z_1,z_2).\eqno(8.10)
$$
The two equations (8.9) are thus mutually consistent. For correlators
with more than two intertwiners one obtains similar equations ([FR]
and below), and for them consistency depends on the fact that $ R$
satisfies the Yang-Baxter relation.

For the other two-point function we have from (5.8), with $z'' = q^{-g-k+uk}z$,
$$\eqalign{ T_pg_{v_0',v_0}(z_1,...,z_N) =
L^{+i}(z_p'')L^{-j}(q^{-g}z_p)
\langle v_0',...\Psi(z_{p-1})L^+_i 
\Psi(z_p)L_j^-\Psi(z_{p+1})...v_0\rangle,\cr}\eqno(8.11)
$$
and, in particular,
$$
T_1g_{v_0',v_0}(z_1,z_2) = L^{+i}(z''_1)L^{-j}(q^{-g}z_1)\langle
v_0',L_i^+
\Psi(z_1)L_j^-\Psi(z_2) v_0\rangle,\eqno(8.12)
$$
$$
T_2g_{v_0',v_0}(z_1,z_2) = L^{+i}(z''_2)L^{-j}(q^{-g}z_2)\langle
 v_0',\Psi(z_1)L^+_i \Psi(z_2)L_j^- v_0\rangle,\eqno(8.13)
$$
and with the help of the Appendix,
$$\eqalign {& T_1g_{v_0',v_0}(z_1,z_2) = q^{A_1}R^{-1}_{12}({z_2\over
z_1})g_{v_0',v_0}(z_1,z_2),\cr & T_2g_{v_0',v_0}(z_1,z_2) =
R_{12}({z_2\over
z_1}q^{-k-g})q^{A_2}g_{v_0',v_0}(z_1,z_2).}\eqno(8.14)
$$

The q-KZ equations for the 3-point functions are
$$
\eqalign{&T_1f(z_1,z_2,z_3) =   R^{-1}_{12}({z_2\over z_1}q^{k+g})
  R^{-1}_{13}({z_3\over z_1}q^{k+g})q^{A_1}f(z_1,z_2,z_3),\cr&
T_2f(z_1,z_2,z_3) = R^{-1}_{23}({z_3\over z_2}q^{k+g})q^{A_2}
R_{12}({z_2\over z_1} )f(z_1,z_2,z_3),\cr & T_3f(z_1,z_2,z_3) =
q^{A_3} R_{13}({z_3\over z_1}) R_{23}({z_3\over z_2}) f(z_1,z_2,z_3),
\cr}\eqno(8.15)
$$
with $A$ as before. Integrability is expressed as a cocycle condition
that is precisely the Yang-Baxter relation for $R$, Eq.(A.6) with $c_2
= 0$.  (The tilde on $\tilde R_{12}$ is redundant.)
\b
\no {\it Remarks 8.} (1) It is interesting to note that
$$
T_1T_2T_3 f(z_1,z_2,z_3) = q^{A_1 + A_2 + A_3}f(z_1,z_2,z_3) =
q^{{1\over 2} (C(\mu) - C(\nu) )}f(z_1,z_2,z_3).\eqno(8.16)
$$
This is what one expects, since the product of any number of
intertwiners should have the universal property; see Remark 3.2, also
(4.31) and (8.10).

(2) The first and the last equations in (8.15) can be written as
follows,
$$\eqalign{ &T_1f(z_1,z_2,z_3) =
q^{-(k+g)d_1}R^{-1}_{1,32}q^{(k+g)d_1}q^{A_1}f(z_1,z_2,z_3),\cr &
T_3f(z_1,z_2,z_3) = q^{A_3} R_{21,3}f(z_1,z_2,z_3).\cr} \eqno(8.17)
$$
Here $R_{i,jk}$ is the action of the universal R-matrix in the
evaluation module via the opposite coproduct, $R_{1,32} = ({\rm id}
\otimes \Delta')R, R_{21,3} = (\Delta' \otimes {\rm id})R$. This too
is an expression of universality; compare the first of (8.17) with the
first of (8.9).
\b
Similarly one finds directly, using the formulas in the Appendix that,
if $g(z_1,z_2,z_3)$ is the alternative 3-point function in (8.1), then
$$
\eqalign{
T_3g(z_1,z_2,z_3)& = R_{23}({z_3\over z_2}q^{-k-g})R_{13}({z_3\over
z_1}q^{-k-g})q^{A_3}g(z_1,z_2,z_3)\cr &=
q^{-(k+g)d_3}R_{12,3}q^{(k+g)d_3}q^{A_3}g(z_1,z_2,z_3).}\eqno(8.18)
$$
The other two formulas cannot be obtained so directly, but the
principle of universality encountered in Remarks 8 tells us that
$$
T_1g(z_1,z_2,z_3) = q^{A_1}R^{-1}_{1,23}g(z_1,z_2,z_3).\eqno(8.19)
$$
Finally, from (8.23),
$$\eqalign{ T_2g(z_1,z_2,z_3) &= T^{-1}_1q^{A_1 +
A_2+A_3}T^{-1}_3g(z_1,z_2,z_3)\cr &=
q^{(k+g)d_1}R_{12}q^{A_2}R^{-1}_{23}q^{-(k+d)d_1}g(z_1,z_2,z_3).}\eqno(8.20)
$$

\b

 Summing up, we have
$$
\eqalign{&T_1g(z_1,z_2,z_3) = q^{A_1}R^{-1}_{1,23}\,g(z_1,z_2,z_3),\cr &
T_2g(z_1,z_2,z_3) =
 q^{-(k+g)d_2}R_{12}q^{(k+g)d_2}q^{A_2}R^{-1}_{23}g(z_1,z_2,z_3),\cr &
 T_3g(z_1,z_2,z_3) = q^{-(k+g)d_3}R_{12,3}q^{(k+g)d_3}q^{A_3}
 g(z_1,z_2,z_3).\cr}\eqno(8.21)
$$

\b
\no{\it Twisting and Covariance.}

Let us evaluate one of the two-point functions of the twisted model,
$$ 
g_\epsilon(z_1,z_2) = \langle
v',\Psi_\epsilon(z_1)\Psi_\epsilon(z_2)\,v\rangle = \langle v',
F_\epsilon^{-1}\Psi
(z_1)F_\epsilon^{-1}\Psi(z_2)\,v\rangle.\eqno(8.22)
$$
On highest weight vectors,
$$
g_\epsilon(z_1,z_2) = \langle
v'_0,\Psi(z_1)F_\epsilon^{-1}(z_2)\Psi(z_2)\,v_0\rangle =
F_\epsilon^{-1,i}(z_2)\langle v'_0,\,\Delta'(F^{-1}_{\epsilon i})
\Psi(z_1)\Psi(z_2)\,v_0\rangle.\eqno(8.23)
$$
Now $\Delta'(f_{-\rho}) = f_{-\rho} \otimes 1 +
q^{\varphi(.,\rho)}\otimes f_{-\rho}$ and so, for Hopf deformations,
when $F_\epsilon$ is a series of the type (6.4),
$$
g_\epsilon(z_1,z_2) = F_\epsilon^{-1}(z_2,z_1)g(z_1,z_2).\eqno(8.24)
$$
An alternative derivation of this result makes direct use of the
cocycle condition. It can be written as follows
$$
F^{-1}_{\epsilon 13}\bigl(({\rm id} \otimes
\Delta_{31})F_\epsilon^{-1}\bigr) = F^{-1}_{\epsilon
21}\bigl((\Delta_{12} \otimes {\rm id})
F_\epsilon^{-1}\bigr).\eqno(8.25)
$$
Applying $v'_0$ we get, because this vector is a highest weight vector, 
$$
 \langle v'_0,({\rm id} \otimes \Delta_{31})F_\epsilon^{-1})... =
\langle v'_0,F^{-1}_{\epsilon 21}...~,\eqno(8.26)
$$
which is just what we need to reduce (8.23) to (8.24).

The transformation formula (8.24) shows that the result (8.14) is not
covariant with respect to twisting, in the following sense.  The
equation satisfied by $g_\epsilon$ is
$$
T_1g_\epsilon(z_1,z_2) = F_\epsilon^{-1}(z_2,q^{-k-g}z_1)
q^{A_1}R^{-1}_{12}({z_2\over z_1})F_\epsilon(z_2,z_1)
g_\epsilon(z_1,z_2);
$$
the right hand side  is very different from the naive analogue of (8.14),
$$
q^{A_1}R^{-1}_{\epsilon 12} ({z_2\over z_1})g_\epsilon(z_1,z_2).
$$
Thus twisting does not preserve the form of the equations satisfied by
 matrix elements of intertwining operators; one cannot simply replace
 $ R $ in these equations by a twisted R-matrix. In fact, it is clear
 that our calculations made use of the specific form of the standard
 R-matrix. The factors $q^A$, in particular, are characteristic of the
 standard R-matrix.  Of course, we do not deny the existence of
 holonomic difference equations that involve R-matrices of a more
 general type. The claim is that the solutions to such equations are
 not, in general, matrix elements of intertwining operators for
 highest weight, quantized Kac-Moody modules.  The elliptic
 correlation functions can be found by solving a ``modified" q-KZ
 equation, but much more simply by the intermediary of the solutions
 of the standard q-KZ equations for the 6-vertex model, as in
 Eq.(8.24).

For three-point functions the effect of twisting is 
$$\eqalign{ g_\epsilon(z_1,z_2,z_3) &= \langle
v_0',F_\epsilon^{-1}\Psi(z_1)F_\epsilon^{-1}\Psi(z_2)F_\epsilon^{-1}
\Psi(z_3)v_0\rangle\cr
& = F_\epsilon^{-1,i}(z_2)F_\epsilon^{-1,j}(z_3)\langle v'_0,
\Psi(z_1)F^{-1}_{\epsilon i}\Psi(z_2)F^{-1}_{\epsilon j}
\Psi(z_3)v_0\rangle\cr & =
F_\epsilon^{-1,i}(z_2)F_\epsilon^{-1,j}(z_3)\langle
v'_0,\Delta_{41}(F^{-1}_{\epsilon
i})\Psi(z_1)\Delta_{42}(F^{-1}_{\epsilon
j})\Psi(z_2)\Psi(z_3)v_0\rangle\cr & =
F_\epsilon^{-1}(z_2,z_1)(F_{\epsilon
j}^{-1})_i(z_2)F_\epsilon^{-1,j}(z_3)\langle
v'_0,\Psi(z_1)(F^{-1}_{\epsilon j})^i\Psi(z_2)\Psi(z_3)v_0\rangle\cr &
= F_\epsilon^{-1}(z_2,z_1)(F_{\epsilon
j}^{-1})_i(z_2)F_\epsilon^{-1,j}(z_3)
\langle v'_0,\Delta_{41}((F^{-1}_{\epsilon j})^i)\Psi(z_1)
\Psi(z_2)\Psi(z_3)v_0\rangle\cr &
= F_\epsilon^{-1}(z_2,z_1)(F_{\epsilon
j}^{-1})_i(z_2)F_\epsilon^{-1,j}(z_3)
\langle v'_0, ((F^{-1}_{\epsilon j})^i)(z_1)\Psi(z_1)\Psi(z_2)
\Psi(z_3)v_0\rangle\cr &
= F_\epsilon^{-1}(z_2,z_1)\Delta_{12}(F_{\epsilon
j}^{-1})(z_1,z_2)F_\epsilon^{-1,j}(z_3)g(z_1,z_2,z_3)
\cr &
= F_\epsilon^{-1}(z_2,z_1)\bigl(({\rm id} \otimes
\Delta)F_\epsilon^{-1}\bigr)(z_3,z_1,z_2)g(z_1,z_2,z_3). 
\cr}\eqno(8.27)
$$ 
Thus we conclude that the twisted correlation functions can be
obtained from the untwisted ones.  The latter are found by solving
equations that are known to be integrable by virtue of the fact that
the standard R-matrix satisfies the Yang-Baxter relation. It is
possible, but redundant and unrewarding, to write down the equations
satisfied by the twisted correlation functions; they are complicated
and uninstructive whether expressed in terms of $R$ or $R_\epsilon$.

\ve
\no{\steptwo 9. Correlation Function for the 8-Vertex Model.}
\b\no
Here we try to understand what, if any, are the qualitative new
features that result from the fact that the elliptic quantum group is
not a Hopf algebra. Technically, the difference is that the reduction
of (8.23) to (8.24) is no longer valid, because the twistor is no
longer of the type (6.4).  Instead of the cocycle condition (6.1) that
gave us (8.25) we now have the modified cocycle condition (7.3), which
yields
$$
F^{-1}_{\epsilon 13}\bigl(({\rm id} \otimes \Delta_{31}
)F_\epsilon^{-1}\bigr) = F^{-1}_{\epsilon 21,3}\bigl((\Delta_{12}
\otimes {\rm id}) F_\epsilon^{-1}\bigr)\eqno(9.1)
$$
and, instead of (8.24), 
$$
g_\epsilon(z_1,z_2) = F^{-1}_{\epsilon 21,3}(z_2,z_1)g(z_1,z_2).\eqno(9.2)
$$
To calculate $F_{\epsilon 21,3}$ see (7.6).  This is the two-point
function for the 8-vertex model. The quasi Hopf nature of the elliptic
quantum group is parameterized by the level $k$ of the highest weight
module and the effect on the two-point function is in the numerical
modification of the matrix $F_\epsilon$ that is indicated by the third
index. Equation (8.27) gets modified in the same manner.

To get an idea of the importance of this effect it is enough to
calculate the modified matrix in the case $N = 2$ with $V$ the
fundamental $\sl(2)$-module.  The result is as follows.  The
trigonometric R-matrix is given for comparison, with the two spaces
interchanged:
$$ 
R^{\rm t} = {A(q,x )\over 1-q^{-2}x}q^\varphi\biggl((1-q^{-2}x )H_+ +
(1-x )H_- + e_1 \otimes e_{-1} + e_0 \otimes e_{-0}\biggr),\eqno(9.3)
$$
$$ 
F_{\epsilon 12,3}^{2m} = {A^{2m}(q,x,\epsilon)\over
1-q^2\alpha\alpha'x}\biggl((1-q^2\alpha\alpha'x)H_+ +
(1-\alpha\alpha'x)H_- - \alpha f_1 \otimes f_{-1} - \alpha'f_0 \otimes
f_{-0}\biggr),\eqno(9.4)
$$
$$
 F_{\epsilon 12,3}^{2m-1} =  { A^{2m-1}(q,x,\epsilon)\over 1-q^2\beta^2 x}
\biggl((1-\beta^2 x)H_+ + (1-q^2\beta^2 x)H_- 
- \beta  f_1 \otimes f_{-0} - \beta f_0 \otimes f_{-1}\biggr).\eqno(9.5) 
$$
Here $ x = z_1/z_2$, 
$$
H_\pm = {\textstyle {1\over 4}}[(1+H)\otimes (1\pm H) + (1- H)\otimes
(1\mp H)]
$$
(in another notation, $H_+ = e_{11}\otimes e_{11} + e_{22}\otimes
e_{22}, H_- = e_{11}\otimes e_{22} + e_{22} \otimes e_{11}$), and
$$ 
\alpha = q^{uk}(\epsilon^2 q^{-k})^m, 
~\alpha' = q^{(1-u)k}(\epsilon^2 q^{-k})^m, ~ \beta^2 =
q^k(\epsilon^2q^{-k})^{2m-1}.
$$
Remember that $k$ denotes the level of the highest weight $\widehat
{\sl(2)}$-module. It enters here because it appears in the extension
$F_{\epsilon 21,3}$ of the twistor in Eq.(9.2). This operator acts in
three spaces, but its action on the highest weight module is limited
to the center. In the level zero case we recover the Hopf twistor.
The calculation that leads to (9.3) is given in detail in the
Appendix, with an explicit formula for the normalizing factor
$A(q,x)$. The matrix factors in (9.4) and (9.5) are obtained in the
same way, and the scalar factors as follows.
\b
 \no {\bf Proposition 9.1.} {\it (a) The normalizing factor in (9.4) is
$$
 A^{2m}(q,x,\epsilon) =  A(1/q ,\alpha\alpha'x).\eqno(9.7)
$$
(b) The normalizing factor in (9.5) is
$$
A^{2m-1}(q,x,\epsilon) = A(1/q,\beta^2x).\eqno(9.8)
$$
}
\b
\no {\it Proof of (a).} Consider the universal R-matrix, 
and the algebra map generated by
$$
e_1 \rightarrow \alpha^{-1} f_1,~~~e_{-1} \rightarrow -\alpha
f_{-1},~~~ e_0 \rightarrow \alpha'^{-1} f_0,~~~ e_{-0} \rightarrow
-\alpha' f_{-0}\eqno(9.9)
$$
in the second space, but $e_i \rightarrow f_i$ in the first
space. This maps the original algebra to another algebra with the $q$
replaced by $q^{-1}$.  Now consider the factorization (2.1) of the
universal R-matrix, $R = q^\varphi T$. After (9.9), the first two
terms in $T^{\rm t}$ agree with the first two terms of
$F_\epsilon^{2m}$.  The recursion relations (2.4) and (7.4) also
agree, after replacing $q$ by $1/q$, and so do the solutions.  Then we
pass to the evaluation representation, setting $f_0 = \hat
z_1f_{-1}$. In the R-matrix (more precisely, in $T^{\rm t}$) we have
set $e_0 = z_1e_{-1}$, which after the substitution (9.9) becomes
$\alpha'f_0 = (z_1/\alpha) f_{-1}$, so that $ z_1 = \alpha\alpha'\hat
z_1$. Under these transformations, including transposition of the two
spaces, the polynomial factor in (9.3) is transformed into that of
(9.4), and the normalizing factor also agree.
\b
\no {\it Proof of (b).} In this case, in the second space let $e_1 
\rightarrow \beta^{-1} f_0,
e_{-1} \rightarrow -\beta f_{-0}, e_0 \rightarrow \beta^{-1}
f_1,\break e_{-0} \rightarrow -\beta f_{-1} $, and in the first space
$e_i
\rightarrow f_i$.
\b

Putting it all together, we have after a simple change of basis
$$
F_{\epsilon 12,3}(z_1,z_2) = A(F_\epsilon)\pmatrix{ a&&& \hat d \cr
&b&\hat c&\cr &\hat c&b&\cr
\hat d  &&& a\cr} ,
\quad x = z_1/z_2,\eqno(9.10)
$$
with 
$$\eqalign{ A(F_\epsilon) = \prod_{m\geq 0}A(q^{-1},q^k \bar
\epsilon^{4m+2}x) &= \prod_{m,n \geq 0} {(1-xq^k \bar
\epsilon^{4m+2}q^{4n})(1-xq^k \bar \epsilon^{4m+2}q^{4n+4})\over
(1-xq^k \bar \epsilon^{4m+2}q^{4n+2})^2} \cr &= {(q^k\bar
\epsilon^2x;q^4,\bar \epsilon^4)_\infty (q^{k+4}\bar
\epsilon^2x;q^4,\bar \epsilon^4)_\infty\over (q^{k+2}\bar
\epsilon^2x;q^4,\bar \epsilon^4)_\infty^2}.\cr}\eqno(9.11)
$$ 
and
$$
a \pm \hat d = \prod_{m\geq 1} { (1 \pm q^{-1+k/2}\sqrt x\,\bar
\epsilon^{2m-1})\over (1\pm q^{1+k/2}\sqrt
x\,\bar\epsilon^{2m-1})},~~~ b \pm \hat c = \prod_{m \geq 1} {(1 \pm
q^{-1+k/2}\sqrt x\,\bar \epsilon^{2m})\over (1\pm q^{1+k/2}\sqrt
x\,\bar\epsilon^{2m})},
$$
with $\bar\epsilon^2 = \epsilon^2q^{-k}$.  Finally, we give the result
of projecting the universal elliptic R-matrix of $\widehat {\sl(2)}$
on the evaluation representation ($k = 0$),
$$
R_\epsilon(z_1,z_2) = \big((F_\epsilon^{\rm
t})^{-1}RF_\epsilon\bigr)(z_1,z_2) =
A_\epsilon(q,x^{-1})\pmatrix{\alpha&&&\delta\cr &\beta&\gamma\cr
&\gamma&\beta\cr d&&&\alpha\cr},
$$
where
$$
\eqalign{
& \alpha+\delta=q{\theta_3(u-\rho,\tau)\over \theta_3(u+\rho,\tau)}, \quad
\alpha-\delta=q{\theta_2(u-\rho,\tau)\over \theta_2(u+\rho,\tau)},\cr
&\beta+\gamma={\theta_1(u-\rho,\tau)\over \theta_1(u+\rho,\tau)},\quad 
\beta-\gamma={\theta(u-\rho,\tau)\over \theta(u+\rho,\tau)},\cr}
$$
with
$$
x = z_1/z_2 = {\rm e}^{4\pi{\rm i}u}, \quad q = {\rm e}^{2\pi{\rm i}\rho},
\quad \epsilon={\rm e}^{\pi{\rm i}\tau},\quad A_\epsilon(q,x^{-1}) =
A(q,x^{-1})A(F_\epsilon)/A(F_\epsilon^{\rm t}).
$$  
In terms of the Jacobian elliptic functions one has
$$
\alpha + \delta : \alpha - \delta : \beta + \gamma : \beta - \gamma 
= {{\rm dn}(2K(u-\rho),k)\over {\rm dn}(2K(u+\rho),k)} : 1 : {1\over
q}{{\rm cn}(2K(u-\rho),k)\over {\rm cn}(2K(u+\rho),k)} : {1\over
q}{{\rm sn}(2K(u-\rho),k)\over {\rm sn}(2K(u+\rho),k)},
$$
where $K, k$ are the real quarter-period and modulus, respectively,
for the nome $\epsilon$:
$$
K={\pi\over 2}\prod_{n\geq 1}\biggl({1+\epsilon^{2n-1}\over
1-\epsilon^{2n-1}}\cdot {1-\epsilon^{2n}\over
1+\epsilon^{2n}}\biggr)^2, \quad k=4\sqrt{\epsilon}\prod_{n\geq
1}\biggl({1+\epsilon^{2n}\over 1+\epsilon^{2n-1}}\biggr)^4.
$$

\bb
\no {\steptwo Acknowledgements.}
\b\no
We thank Olivier Babelon, Benjamin Enriquez, Moshe Flato, Tetsuji Miwa
and Nikolai Reshe\-ti\-khin for advice. We thank Moshe Flato for an
incisive and constructive criticism of the original manuscript.
A.G. thanks the Fundaci\'on Del Amo for financial support and the
Department of Physics of UCLA for hospitality.

\bb
\ce {\steptwo Appendix.}
\b
\no {\it Solving the recursion relations.}
 
We shall solve the recursion relation (2.4) in the fundamental
evaluation representation of $\widehat{\sl(2)}$. Here we set
$$
e_1 = \kappa\pmatrix{ 0&1\cr0&0\cr}, ~~~e_{-1} =
\kappa\pmatrix{0&0\cr1&0\cr},~~~
\varphi = {1\over 2}H \otimes H,~~~ H = \pmatrix{1&0\cr 0&-1}.\eqno(\rm A.1) 
$$
The commutation relations hold with $\kappa^2 = q - q^{-1}$. The
factor $T$ in $R = q^\varphi T$ has the form
$$
T = \pmatrix{a\cr&b&cx\cr &c&b\cr&&&a\cr},~~~x = z_1/z_2,
$$ 
and (2.4) is equivalent to
$$
[T,1 \otimes e_{-\gamma}] = \bigl(e_{-\gamma} \otimes
q^{\varphi(\gamma,.)}\bigr)T - T\bigl(e_{-\gamma} \otimes
q^{-\varphi(.,\gamma)}\bigr),~~~\gamma = 1,0,
$$  
with $\varphi(1,.) = \varphi(.,1) = H, \varphi(0,.) = \varphi(.,0) = -H$. 
Taking $\gamma = 1$ we get two relations,
$$
q(a-b) = c = (aq - b/q)/x,\eqno(\rm A.2)
$$
and taking $\gamma = 0$ the same two relations. Hence 
$$ 
R(q,x^{-1}) = {A(q,x^{-1})\over
1-q^{-2}x^{-1}}q^\varphi\bigl((1-q^{-2}x^{-1})H_+ + (1-x^{-1})H_- +
e_{-1} \otimes e_1 + e_{-0} \otimes e_0\bigr).
$$
The matrices in (9.4) and (9.5) are found in the same way. In the
special case of $\widehat {\sl(2)}$ the Cartan factors in (7.5) are,
for $m \geq 1$,
$$
Q(2m,1) = q^{(u-m)c},\quad Q(2m,0) = q^{(1-u-m)c},\quad Q(2m-1,1) =
Q(2m-1,0) = q^{(1-m)c}.
$$

In the structure, $R$ is determined uniquely by the recursion
relations and the initial conditions, but in the evaluation
representation the normalizing factor $A(q,x^{-1})$ remains
undetermined.  Fortunately Levendorskii, Soibelman and Stukopin [LS],
starting from an equivalent expression for the standard, universal
R-matrix for $\widehat {\sl(2)}$ obtain the following result,
$$
A(q,x) = \exp\biggl(\sum_{k \geq 1}{1\over k}{q^k - q^{-k} \over q^k +
q^{-k}}x^k\biggr).\eqno(\rm A.3)
$$
The sum converges for $|q| \neq 1, |x| <1$ and the formula can be
manipulated to yield
$$
A(q,x) = \cases{
\displaystyle { (xq^2;q^4)_\infty^2\over 
(x;q^4)_\infty(xq^4;q^4)_\infty}, & $|q| < 1$,\cr
\displaystyle{(x;q^{-4})_\infty(xq^{-4};q^{-4})_\infty\over 
(xq^{-2};q^{-4})_\infty^2},& $|q| >1$.\cr}\eqno(\rm A.4)
$$
Hence 
$$
A(q,x)A(q^{-1},x) = 1, ~~~|q| \neq 1.\eqno(\rm A.5)
$$
This is also clear from (A.3).

The inverse of $R$ can also be represented as a series, similar to (2.1),
$$
R^{-1} = q^{-\varphi}\hat T,~~~\hat T = 1 + \sum_\alpha \hat
e_{-\alpha} \otimes \hat e_\alpha + ...\, ,
$$
with
$$
\hat e_\alpha := q^{-\varphi(\alpha,.)}e_\alpha,~~~ \hat e_{-\alpha} = 
- e_{-\alpha}q^{\varphi(.,\alpha)}.
$$
The commutation relations for the $\hat e$'s agree with those of the
$e$'s, and the recursion relations for $\hat T$ agrees with that of
$T$, all up to the sign of $\varphi$. (We get a recursion relation for
$\hat T$ from the fact that $R^{-1}$ also satisfies the Yang-Baxter
relation.) Consequently, in the structure,
$$
R(\varphi,e)^{-1} = R(-\varphi,\hat e),
$$
and in any evaluation representation,
$$
R(q,x)^{-1} = R(q^{-1},x).
$$ 
These results are quite general and imply, in particular, Eq.(A.5).
\b
\no {\it Reduced formulas.} 

We list here the formulas that are obtained from the Yang-Baxter
relation
$$
R_{12}R_{13}R_{23} = R_{23}R_{13}R_{12}
$$
and the quasi triangular conditions
$$
({\rm id} \otimes \Delta)R = R_{12}R_{13},\quad (\Delta \otimes {\rm
id})R = R_{23}R_{13}
$$
when the $c,d$ factors are removed as in
$$
R = q^{uc\,\otimes \,d + (1-u)d \,\otimes \,c}\tilde R,
$$
namely
$$
\tilde R_{12}(q^{-uc_2d_3}\tilde R_{13}q^{uc_2d_3})\tilde R_{23} = 
\tilde R_{23}(q^{-(1-u)d_1c_2}\tilde R_{13}q^{(1-u)d_1c_2})\tilde 
R_{12}\eqno(\rm A.6) 
$$
and
$$
({\rm id} \otimes \Delta)\tilde R = (q^{-(1-u)d_1c_3}\tilde
R_{12}q^{(1-u)d_1c_3})\tilde R_{13},\quad (\Delta \otimes {\rm
id})\tilde R = (q^{-uc_1d_3}\tilde R_{23}q^{uc_1d_3})\tilde
R_{13}.\eqno(\rm A.7)
$$
These last two relations give us what we need to reduce (8.8), namely
$$
 L^-(z_1)_{13}\Phi(z_2) = 
\bigl(\tilde R_{12}({z_2\over z_1}q^{k-uk})\bigr)^{-1}\bigl(L^{-i}(z_1)\otimes
1\bigr)\Phi(z_2)L^-_i,\eqno(\rm A.8)
$$
$$ L^{+i}(z_2)\Phi(z_1)L_i^+ = 
L^+(z_2)\tilde R_{12}({z_2\over z_1})\Phi(z_1). \eqno(\rm A.9) 
$$
For the other intertwiner, there is an analogue of (A.9),
$$
L^{+i}(z_2)\Psi(z_1)L^+_i = \tilde R_{12}({z_2\over
z_1}q^{-uk})L^+(z_2)\Psi(z_1),\eqno(\rm A.10)
$$
but we could not find an analogue of (A.8).  To obtain (8.14) we used
the method that was explained for the derivation of (8.21).

\bb
 
\no{{\bf \steptwo References.}}
\b 

\item {[BBB]} O. Babelon and D. Bernard, A Quasi-Hopf interpretation of  
 quantum 3-$j$ and 6-$j$ symbols and difference equations,
q-alg/9511019.

\item {[Ba]} R.J. Baxter, Partition Function of the Eight-Vertex Model,  Ann.
Phys.  {\bf 70}  (1972) 193-228.

\item {[BK]} R.J. Baxter and S.B. Kelland, J. Phys. C: Solid State Phys. 
{\bf 7} (1974) L403-6.

\item {[Be]} A.A. Belavin, Dynamical Symmetry of Integrable Systems, 
Nucl. Phys.
{\bf 180} (1981) 198-200.

\item {[BD]} A.A. Belavin and V.G. Drinfeld, Triangle Equation and Simple Lie
Algebras,   Sov. Sci. Rev. Math. {\bf 4} (1984) 93-165.

\item {[Ber]} D. Bernard, On the WZW model on the torus, Nucl. Phys. 
{\bf B303} (1988) 77-174. 

\item {[D1]}V.G. Drinfeld, Quantum Groups, in Proceedings, 
International Congress of Mathematicians, Berkeley, A.M. Gleason,
 ed. (A.  M.  S., Providence 1987).

\item {[D2]} V.G. Drinfeld, Quasi Hopf Algebras, Leningrad Math. J. {\bf 1} 
(1990) 1419-1457.

\item {[ER]} B. Enriquez and   Rubtsov, Quasi-Hopf algebras associated 
with sl(2) and complex curves,
\break
q-alg/9703018. 

\item {[Fe]} G. Felder, Elliptic Quantum Groups, hep-th/9412207.

 \item {[FR]} I.B. Frenkel and N.Yu. Reshetikhin, Quantum Affine
 Algebras and Holonomic Difference Equations, Commun. Math. Phys.
 {\bf 146} (1992) 1-60.

\item {[FRS]} I.B. Frenkel, N.Yu. Reshetikhin and M. Semenov-Tian-Shansky,
 Drinfeld-Sokolov reduction for difference operators and deformations
     of W-algebras I. The case of Virasoro algebra, q-alg/9704011.

\item {[Fr1]} C. Fr\o nsdal, Generalization and Deformations of 
Quantum Groups, to appear in RIMS Publications. (q-alg/9606020)

\item {[Fr2]} C. Fr\o nsdal, Quasi Hopf Deformation of Quantum Groups, 
to appear  in Letters in Mathematical Physics, q-alg/9611028. 

\item {[JM]} M. Jimbo and T. Miwa, Algebraic Analysis of Solvable Lattice 
Models, Regional Conference Series in Mathematics, (1995) Number 85.

\item {[JMN]}   M. Jimbo, T. Miwa and A. Nakayashiki, Difference equations 
for the correlation functions of the eight-vertex model,
J. Phys. A: Math. Gen.  {\bf  206 } (1993) 2199-2209.

\item {[KZ]} V.G. Knizhnik and A.B. Zamolodchikov, Current algebra and 
Wess-Zumino model in two dimensions, Nucl. Phys. B {\bf 247} (1984) 83-103.

\item {[LS]} S. Levendorskii, Y. Soibelman and V. Stukopin, The Quantum Weyl 
Group and the Universal Quantum R-Matrix for Affine Lie Algebra
$A_1^{(1)}$, Lett. Math. Phys. {\bf 27} (1993) 253-264.

\item {[R]} N.Yu. Reshetikhin, Multiparameter Quantum Groups and Twisted 
Quasitriangular Hopf Algebras, Lett. Math. Phys. {\bf 20} (1990)
331-336.
\ve

\end